\documentclass{article}
\usepackage{amssymb}
\usepackage{amsmath}
\usepackage{amsfonts}
\usepackage{graphicx}

\newtheorem{theorem}{Theorem}[section]
\newtheorem{lemma}[theorem]{Lemma}
\newtheorem{corollary}[theorem]{Corollary}

\newtheorem{remark}[theorem]{Remark}

\newcommand{\qed}{\hfill \rule{2.3mm}{2.3mm}}

\newcommand{\om}{\omega}
\newcommand{\al}{\alpha}

\newcommand{\la}{\lambda}

\newcommand{\ve}{\varepsilon}
\newcommand{\vp}{\varphi}

\newcommand*{\esssup}{\operatorname*{ess\phantom{|}\!sup}}
\newcommand*{\essinf}{\operatorname*{ess\phantom{|}\!inf}}

\newcommand{\bee}{\begin{equation}}
\newcommand{\ee}{\end{equation}}
\newcommand{\uu}{\bar{u}_\ve}

\newcommand{\R}{\mathbb{R}}
\newcommand{\N}{\mathbb{N}}
\newcommand{\M}{\mathbb{M}}
\newcommand{\Z}{\mathbb{Z}}

\newcommand{\x}{x/\varepsilon}

\newcommand{\LL}{{\cal L}}
\newcommand{\reff}[1]{(\ref{#1})} 
  
\begin{document}

\title{\bf
A common approach to singular perturbation and homogenization III: 
Nonlinear periodic homogenization with localized defects
\footnote{\bf This work is dedicated to the memory of Nikolai N. Nefedov, the coauthor of parts I and II of the present paper, a great mathematician and a wonderful friend.
He passed away in October 2024.}
}

\newcounter{thesame}
\setcounter{thesame}{1}
\author{
Lutz Recke \thanks{Humboldt University of Berlin, Institute of Mathematics, Rudower Chaussee 25, 12489 Berlin, Germany.
		{\small   E-mail:
			{\tt lutz.recke@hu-berlin.de}}
}}
\date{}

\maketitle

\begin{abstract}
\noindent
We consider  periodic homogenization with localized defects for semilinear 
elliptic equations and systems of the type
$$
\nabla\cdot\Big(\Big(A(\x)+B(\x)\Big)\nabla u(x)+c(x,u(x)\Big)=d(x,u(x)) \mbox{ in }  \Omega
$$
with Dirichlet boundary conditions.
For small $\ve>0$ we show existence of weak solutions $u=u_\ve$ as well as their local uniqueness for $\|u-u_0\|_\infty \approx 0$,
where $u_0$ is a given
non-degenerate weak solution to the homogenized problem. Moreover, we prove that $\|u_\ve-u_0\|_\infty\to 0$ for $\ve \to 0$, and we estimate the corresponding rate of convergence.

Our assumptions are, roughly speaking, as follows: $\Omega$ is a bounded Lipschitz domain, $A$, $B$, $c(\cdot,u)$ and $d(\cdot,u)$ are bounded and measurable, 
$c(x,\cdot)$ and $d(x,\cdot)$ are $C^1$-smooth, $A$ is periodic, and $B$ is a localized defect.
Neither global uniqueness is supposed nor growth restriction for $c(x,\cdot)$ or $d(x,\cdot)$. 

The main tool of the proofs is an abstract result of implicit function theorem type
which permits a common approach to nonlinear singular perturbation and homogenization.
\end{abstract}



\section{Introduction}
\setcounter{equation}{0}
\setcounter{theorem}{0}

This paper concerns periodic homogenization with localized defects in the sense of \cite{Blanc}. More exactly, we consider Dirichlet
problems for $2D$
semilinear second-order elliptic systems of the type
\bee
\label{BVP1}
\left.
\begin{array}{l}
\partial_{x_i}\Big(\big(a^{\al \beta}_{ij}(\x)
+b^{\al \beta}_{ij}(\x)\big)
\partial_{x_j}u^\beta(x)
+c_i^\al(x,u(x))\Big)=d^\al(x,u(x)) \mbox{ in }  \Omega,\\
u^\al(x)=0\mbox{ on } \partial\Omega,
\end{array}
\right\}
\al=1,\ldots,n
\ee
as well as for $ND$ semilinear second-order elliptic equations of the type
\bee
\label{BVP2}
\left.
\begin{array}{l}
\partial_{x_i}\Big(\big(a_{ij}(\x)+b_{ij}(\x)\big)\partial_{x_j}u(x)+c_i(x,u(x))\Big)= d(x,u(x)) \mbox{ in } \Omega,\\
u(x)=0
\mbox{ on } \partial\Omega.
\end{array}
\right\}
\ee
Here $\ve>0$ is the small homogenization parameter.
In \reff{BVP1} repeated indices are to be summed over $i,j,\ldots=1,2$ and $\al,\beta,\ldots=1,\ldots,n$, and in \reff{BVP2}
repeated indices are to be summed over $i,j,\ldots=1,\ldots,N$.
Concerning \reff{BVP1} we suppose that
\begin{eqnarray}
\label{Omass1}
&&\mbox{$\Omega$ is a bounded Lipschitz domain in $\R^2$,}\\
\label{aass1}
&&a_{ij}^{\al \beta} \in L^\infty(\R^2)
\mbox{ are $\Z^2$-periodic, } b^{\al \beta}_{ij} \in L^\infty(\R^2)\cap L^1(\R^2),\\
\label{mona1}
&&\essinf\left\{a^{\al \beta}_{ij}(y)v^\al_i v^\beta_j:\; y\in \R^2,\,v \in \R^{2n},\; v^\al_iv^\al_i=1\right\}>0,\\
\label{monb1}
&&\essinf\left\{\big(a^{\al \beta}_{ij}(y)+b^{\al \beta}_{ij}(y)\big)v^\al_i v^\beta_j:\; y\in \R^2,\,v \in \R^{2n},\; v^\al_iv^\al_i=1
\right\}>0,\\
\label{diffass1}
&&u\in \R^n \mapsto (c_i^\al(\cdot,u),
d^\al(\cdot,u))
\in L^\infty(\Omega)^2 \mbox{ are $C^1$-smooth}.
\end{eqnarray}
And similarly, concerning \reff{BVP2} we suppose that 
\begin{eqnarray}
\label{Omass2}
&&\Omega \mbox{ is a bounded Lipschitz domain in } \R^N \mbox{ with } N \in \N,\\
\label{aass2}
&&a_{ij} \in L^\infty(\R^N) \mbox{ are $\Z^N$-periodic, } b_{ij} \in L^\infty(\R^N) \mbox{ have compact supports,}\\
\label{monb2}
&&\essinf\left\{\big(a_{ij}(y)+b_{ij}(y)\big)v_i v_j:\; y,v \in \R^N,\; v_iv_i=1\right\}>0,\\
\label{diffass2}
&&u \in \R \mapsto (c_i(\cdot,u),
d(\cdot,u)) \in  L^\infty(\Omega)^2 \mbox{ are $C^1$-smooth.}
\end{eqnarray}

We are going to prove existence and local uniqueness of
weak solutions $u$ to \reff{BVP1} and to 
\reff{BVP2}, respectively,
with $\ve \approx 0$ and $\|u-u_0\|_\infty\approx 0$, where $u_0$ is a given non-degenerate solution to the corresponding homogenized boundary value problem.
The homogenized boundary value problem corresponding to \reff{BVP1} is
\bee
\label{hombvp1}
\left.
\begin{array}{l}
\partial_{x_j}\Big(\hat{a}^{\al \beta}_{ij}
\partial_{x_j}u^\beta(x)
+c_i^\al(x,u(x))\Big)
=d^\al(x,u(x)) \mbox{ in } \Omega,\\
u^\al(x)=0 \mbox{ on } \partial\Omega
\end{array}
\right\}
\al=1,\ldots,n
\ee
with
\bee
\label{hatAdef1}
\hat{a}^{\al\beta}_{ij}:=
\int_{(0,1)^2}\left(a^{\al \beta}_{ij}(y)
+a^{\al \gamma}_{ik}(y)
\partial_{y_k}v^{\gamma \beta}_j(y)\right)dy,
\ee
where the $2n$ correctors $v^{\beta}_j=\left(v^{1 \beta}_j,\ldots,v^{n \beta}_j\right)
\in W^{1,2}_{\rm loc}(\R^2;\R^n)$, $\beta=1,\ldots,n$, $j=1,2$,
are defined by the $2n$ cell problems
\bee
\label{cell1}
\left.
\begin{array}{l}
\partial_{y_i}
\left(a^{\al \beta}_{ij}(y)+
a^{\al \gamma}_{ik}(y)
\partial_{y_k}v^{\gamma \beta}_j(y)\right)
=0\mbox{ for } y\in \R^2,\\
v^{\al \beta}_j(\cdot+z)=v^{\al \beta}_j \mbox{ for } z \in \Z^2,\;
\displaystyle\int_{(0,1)^2} v^{\al \beta}_j(y)dy=0,
\end{array}
\right\}
\al=1,\ldots,n.
\ee
And similarly,
the homogenized boundary value problem corresponding to \reff{BVP2} is
\bee
\label{hombvp2}
\left.
\begin{array}{l}
\partial_{x_j}\Big(\hat{a}_{ij}
\partial_{x_j}u(x)
+c_i(x,u(x))\Big)
=d(x,u(x)) \mbox{ in } \Omega,\\
u(x)=0 \mbox{ on } \partial\Omega
\end{array}
\right\}
\ee
with
\bee
\label{hatAdef2}
\hat{a}_{ij}:=
\int_{(0,1)^N}\Big(a_{ij}(y)
+a_{ik}(y)
\partial_{y_k}v_j(y)\Big)dy,
\ee
where the $N$ correctors $v_j
\in W^{1,2}_{\rm loc}(\R^N)$,  $j=1,\ldots,N$,
are defined by the $N$ cell problems
\bee
\label{cell2}
\left.
\begin{array}{l}
\partial_{y_i}
\Big(a_{ij}(y)+
a_{ik}(y)
\partial_{y_k}v_j(y)\Big)
=0\mbox{ for } y\in \R^N,\\
v_j(\cdot+z)=v_j \mbox{ for } z \in \Z^N,\;
\displaystyle
\int_{(0,1)^N} v_j(y)dy=0.
\end{array}
\right\}
\ee

As usual, a function $u \in W_0^{1,2}(\Omega;\R^n)\cap L^\infty(\Omega;\R^n)$ is called weak solution to the boundary value problem \reff{BVP1} if it satisfies the variational equation
$$
\int_\Omega
\left(\Big(\big(a^{\al \beta}_{ij}(\x)
+b^{\al \beta}_{ij}(\x)\big)
\partial_{x_j}u^\beta(x)+c_i^\al(x,u(x))\Big)
\partial_{x_i}\vp^\al(x)+d^\al(x,u(x))\vp^\al(x)\right)dx=0
$$
for all 
$\vp \in W_0^{1,2}(\Omega;\R^n)$,
and similarly for \reff{BVP2},
\reff{hombvp1} and \reff{hombvp2}. 

Besides the homogenized boundary value problems \reff{hombvp1} and \reff{hombvp2} we consider also their linearizations in a function $u_0$, i.e.
\bee
\label{linhombvp1}
\left.
\begin{array}{l}
\partial_{x_j}\Big(\hat{a}^{\al \beta}_{ij}
\partial_{x_j}u^\beta(x)
+\partial_{u_\gamma}c_i^\al(x,u_0(x))u^\gamma(x)\Big)
=\partial_{u_\gamma}d^\al(x,u_0(x))u^\gamma(x) \mbox{ in } \Omega,\\
u^\al(x)=0 \mbox{ on } \partial\Omega
\end{array}
\right\}
\al=1,\ldots,n
\ee
and
\bee
\label{linhombvp2}
\left.
\begin{array}{l}
\partial_{x_j}\Big(\hat{a}_{ij}
\partial_{x_j}u(x)
+\partial_{u}c_i(x,u_0(x))u(x)\Big)
=\partial_{u}d(x,u_0(x))u(x) \mbox{ in } \Omega,\\
u(x)=0
\mbox{ on } \partial\Omega.
\end{array}
\right\}
\ee
As usual, a function $u \in W_0^{1,2}(\Omega;\R^n)$ is called weak solution to the boundary value problem \reff{linhombvp1} if it satisfies the variational equation
$$
\int_\Omega
\left(\Big(\hat a^{\al \beta}_{ij}
\partial_{x_j}u^\beta(x)+\partial_{u^\gamma}c_i^\al(x,u_0(x))
u^\gamma(x)\Big)
\partial_{x_i}\vp^\al(x)+
\partial_{u^\gamma}d^\al(x,u_0(x))u^\gamma(x)\vp^\al(x)\right)dx=0
$$
for all 
$\vp \in W_0^{1,2}(\Omega;\R^n)$,
and similarly for \reff{linhombvp2}.

And finally, as usual, we denote by $\|\cdot\|_\infty$ the norms in the Lebesgue spaces $L^\infty(\Omega)$ and $L^\infty(\Omega;\R^n)$, i.e.
\bee
\label{normdef}
\|u\|_\infty:=
\left\{
\begin{array}{l}
\displaystyle\esssup_{x \in \Omega}|u(x)|
\mbox{ for } u \in L^\infty(\Omega),\\ 
\displaystyle\sum_{\al=1}^n \esssup_{x \in \Omega}|u^\al(x)|
\mbox{ for } u \in L^\infty(\Omega;\R^n).
\end{array}
\right.
\ee\\

Our result concerning \reff{BVP1} is the following
\begin{theorem} 
\label{main1}
Suppose \reff{Omass1}-\reff{diffass1}, and
let $u=u_0$ be a weak solution to \reff{hombvp1} such that
\reff{linhombvp1}
does not have weak solutions $u\not=0$.
Then the following is true:

(i) There exist $\ve_0>0$ and  $\delta>0$
such that for all $\ve \in (0,\ve_0]$ there exists exactly one 
weak solution $u=u_\ve$ 
to \reff{BVP1} with $\|u-u_0\|_\infty \le \delta$. Moreover, $\|u_\ve-u_0\|_\infty\to 0$
for $\ve \to 0$.

(ii) Suppose that $u_0 \in W^{2,p_0}(\Omega;\R^n)$
with certain $p_0>2$. Then there exists $\la>0$ such that $\|u_\ve-u_0\|_\infty=O(\ve^\la)$
for $\ve \to 0$. If additionally $b_{ij}^{\al \beta}=0$ for all $\al,\beta=1,\ldots,n$ and $i,j=1,2$, then for all $\la \in (0,1/2)$ we have that $\|u_\ve-u_0\|_\infty=O(\ve^\la)$
for $\ve \to 0$.
\end{theorem}

And similarly, our result concerning \reff{BVP2} is 
\begin{theorem} 
\label{main2}
Suppose \reff{Omass2}-\reff{diffass2}, and
let $u=u_0$ be a weak solution to \reff{hombvp2} such that
\reff{linhombvp2}
does not have weak solutions $u\not=0$.
Further, suppose that $u_0 \in W^{2,p_0}(\Omega)$
with certain $p_0>N$.
Then the following is true:

(i) There exist $\ve_0>0$ and  $\delta>0$
such that for all $\ve \in (0,\ve_0]$ there exists exactly one 
weak solution $u=u_\ve$ 
to \reff{BVP2} with $\|u-u_0\|_\infty \le \delta$. Moreover,
there exists $\la>0$ such that $\|u_\ve-u_0\|_\infty=O(\ve^\la)$
for $\ve \to 0$. 

(ii) If  $b_{ij}=0$ for all  $i,j=1,\ldots,N$, then 
for all $\la \in (0,1/N)$ we have that $\|u_\ve-u_0\|_\infty=O(\ve^\la)$
for $\ve \to 0$.
\end{theorem}

The main goal of the present paper is to prove Theorems \ref{main1} amd \ref{main2}
by means of an abstract approach which leads to existence, local uniqueness and estimates of solutions to parameter depending equations, i.e. by means of a result of implicit function theorem type (cf. Theorem \ref{ift} and Corollary \ref{cor} below). Contrary to the classical implicit function theorem, in this approach the equations are allowed to depend on the parameter in a rather singular way, like in singularly perturbed problems or in homogenization problems. In other words: The classical implicit function theorem cannot be applied to the boundary value problems \reff{BVP1} and \reff{BVP2} because the parameter~$\ve$ enters into \reff{BVP1} and \reff{BVP2} in a singular way, but Theorem \ref{ift} and Corollary \ref{cor} can be applied (cf. also Remark \ref{remsing} below).

Remark that not the non-smoothness of the data in \reff{BVP1} and \reff{BVP2}  is the reason why the classical implicit function theorem cannot be applied. It is well-known how to apply the classical implicit function theorem to "regularly perturbed" elliptic and parabolic PDEs with non-smooth data (cf. \cite{GREvol,GrR,PRS,Recke1995}). The non-smoothness of the data in \reff{BVP1} and \reff{BVP2} does not influence the abstract approach to prove Theorems \ref{main1} and \ref{main2}, it influences only the choice of the function spaces which have to be used.

\begin{remark}
It is well-known (as a consequence of the
assumptions \reff{aass1} and \reff{mona1} of Theorem~\ref{main1} and of the
Lax-Milgram lemma, cf. \cite[Section 2.2 and Lemma 2.2.4]{Shen}) that  the cell problem \reff{cell1}
is uniquely weakly solvable  
and that the homogenized diffusion coefficients $\hat{a}^{\al \beta}_{ij}$ 
satisfy the coercivity condition
\reff{mona1} also. 
Therefore,
the rather implicit assumption of Theorem \ref{main1}, that there do not exist nonzero weak solutions 
to \reff{linhombvp1}, is satisfied, for example,  if 
$$
\essinf\left\{\partial_{u^\beta}d^{\al}(x,u_0(x))v^\al v^\beta:\; x \in \Omega,\; v \in \R^{n},\; v^\al v^\al=1\right\}\ge 0
$$
and if $\|\partial_{u^\beta}c_i^{\al}(\cdot,u_0(\cdot))\|_\infty$ are sufficiently small.

And similarly for the cell problem \reff{cell2} and the assumption of Theorem \ref{main2} that there do not exist nonzero weak solutions 
to \reff{linhombvp2}. Remark that in Theorem \ref{main2} we do not assume that
\bee
\label{mona2}
\essinf\left\{a_{ij}(y)v_i v_j:\; y,v \in \R^N,\; v_iv_i=1\right\}>0,
\ee
because \reff{mona2} follows from the assumptions \reff{aass2} and \reff{monb2} of Theorem~\ref{main2}. Indeed, we have
\begin{eqnarray*}
&&\essinf\left\{\big(a_{ij}(y)
+b_{ij}(y)\big)
v_i v_j:\; y,v\in \R^N,\; v_iv_i=1\right\}\\
&&\le\essinf\left\{\big(a_{ij}(y)
+b_{ij}(y+z)\big)
v_i v_j:\; y\in (0,1)^N,\,v \in \R^{N},\; v_iv_i=1\right\}
\mbox{ for any } z \in \Z^N.
\end{eqnarray*}
Because the functions $b_{ij}$ have compact supports, it follows, if  $z_i z_i$ is large, that
\begin{eqnarray*}
&&\essinf\left\{\big(a_{ij}(y)
+b_{ij}(y)\big)
v_i v_j:\; y,v\in \R^N,\; v_iv_i=1\right\}\\
&&\le\essinf\left\{a_{ij}(y)
v_i v_j:\; y,v\in \R^N,\; v_iv_i=1\right\}.
\end{eqnarray*}
\end{remark}

\begin{remark}
Theorem \ref{main1} is true also for mixed boundary conditions of the type
$$
\left.
\begin{array}{l}
\Big(\big(a^{\al \beta}_{ij}(\x)+b^{\al \beta}_{ij}(\x)\big)\partial_{x_j}u^\al(x)+c^\al_i(x,u(x))\Big)\nu_i(x)=
d^\al_0(x,u(x))
\mbox{ on } \Gamma,\\
u^\al(x)=0  \mbox{ on } \partial \Omega \setminus \Gamma,
\end{array}
\right\}
\al=1,\ldots,n
$$
where $\Gamma$ is a  subset of $\partial \Omega$ such that its relative to $\partial \Omega$ boundary consists of finitely many points, and if the maps
$u\in \R^n\mapsto d_0^\al(\cdot,u) \in L^\infty(\Gamma)$ are $C^1$-smooth (cf. also
Remark \ref{boundaryterms} below).
And similarly for Theorem \ref{main2}, where the
relative to $\partial \Omega$ boundary of
$\Gamma$ should be a Lipschitz hypersurface in $\partial \Omega$.
The reason for that is that
the maximal regularity results Theorem \ref{maxreg1} and 
Theorem \ref{maxreg2}
below are true for those boundary conditions also.

In \cite{II} is presented a result
of the type of Theorem \ref{main1} 
for 2D semilinear elliptic systems with various (including mixed) boundary conditions, but without localized defects.

Also, Theorems \ref{main1} and \ref{main2} are true not only for Lipschitz domains
$\Omega$, 
but also for domains which are regular in the sense of \cite[Definition 2]{G}.
The class of regular domains is slightly larger than the class of Lipschitz domains. In particular, bi-Lipschitz transformations of regular domains are regular domains again.
\end{remark}

\begin{remark}
We do not believe that in the case of
$N>2$ space dimensions
Theorem \ref{main1} 
is true for   general elliptic systems (i.e. for systems with  system dimension $n>1$ and with large cross diffusion), no matter if localized defects are present or not.
The reasons for this guess are the well-known examples of unbounded weak solutions to linear elliptic systems 
with  $L^\infty$-coefficients and
with $n>1$, 
$N>2$ and with large cross diffusion (cf., e.g. \cite[Section 8.7]{BenF},
\cite[Section 12.2]{Chen},
\cite[Section 6.2]{Giusti}).

On the other hand, we believe that 
Theorem \ref{main2} can be generalized to the case of elliptic systems which are close to be triangular, i.e. to elliptic systems with diffusion coefficients $a^{\al \beta}_{ij}$ with $\|a^{\al \beta}_{ij}\|_\infty\approx 0$ for $\al>\beta$,
because 
Theorem \ref{maxreg2} below is true for those systems also.
\end{remark}

\begin{remark}
There is a certain asymmetry between Theorems \ref{main1} and \ref{main2}: In Theorem \ref{main1}(i)
existence and local uniqueness is stated without any additional regularity assumption concerning $u_0$, but in Theorem \ref{main2}(i)
existence and local uniqueness is stated under the additional regularity assumption $u_0 \in W^{2,p_0}(\Omega)$ with certain $p_0>N$, only. We do not know if the assertion of Theorem~\ref{main2}(i) remains true if this additional regularity assumption is canceled.

The reason for this asymmetry in the present paper is the used technique:
Theorem~\ref{main1} is proved by working in Lebesgue spaces, in which the set of smooth functions is dense. Therefore the property \reff{Sconv} below can be used.
But Theorem~\ref{main2} is proved by working in Morrey spaces, in which the set of smooth functions is not dense, and no replacement of \reff{Sconv} seems to be available.

The reason why we do not work in Lebesgue and Sobolev spaces in the proof of Theorem~\ref{main2} can be explained easily: The maximal regularity result Theorem \ref{maxreg1} yields maximal regularity in Sobolev spaces $W^{1,p}$ with $p>2$, but $p\approx 2$ only, and these Sobolev spaces are not embedded into $L^\infty$ for space dimension $N>2$ and, hence, on those function spaces the nonlinearities are not well-defined or not $C^1$-smooth, in general.
\end{remark}

\begin{remark}
If we consider the boundary value problem \reff{BVP1} with smooth data and vanishing localized defects, then Theorem \ref{main1}(ii) claims that
$$
\|u_\ve-u_0\|_\infty=O(\ve^\la)  
\mbox{ for } \ve \to 0
\mbox{ for any } \la \in (0,1/2).
$$
Even in the linear case,
i.e. for periodic homogenization of general linear $2D$ elliptic systems with smooth data and Dirichlet boundary conditions,
we do not know if this convergence rate can be improved. 

If we consider the boundary value problem \reff{BVP2} with smooth data and vanishing localized defects, then Theorem \ref{main2}(ii) claims that
$$
\|u_\ve-u_0\|_\infty=O(\ve^\la)  
\mbox{ for } \ve \to 0
\mbox{ for any } \la \in (0,1/N).
$$
But for periodic homogenization of linear $ND$ scalar elliptic equations on smooth domains and with appropriate right-hand sides and Dirichlet boundary conditions it can be shown by means of the maximum principle that
\bee
\label{rate2}
\|u_\ve-u_0\|_\infty=O(\ve)  
\mbox{ for } \ve \to 0
\ee
(cf. \cite[Section 2.4]{Ben}, \cite[Theorem 3.4]{Kenig}).
We do not know if the rate \reff{rate2} is true in general under the assumptions of Theorem \ref{main2}(ii) also.
\end{remark}

\begin{remark}
The cell problems \reff{cell1} and \reff{cell2} and, hence, the corresponding correctors and homogenized diffusion tensors do not depend on the localized defects $b_{ij}^{\al \beta}$ and $b_{ij}$, respectively. Therefore also the solutions $u_0$ to the homogenized problems, which are used in Theorems \ref{main1} and \ref{main2}, do not depend on the localized defects. 

Remark that, if one works with correctors, which depend on the localized defects, and if one assumes higher regularity of the data of  \reff{BVP1} and  \reff{BVP2}, then more precise estimates of the solution family $u_\ve$ can be obtained (cf. \cite{BJL} for the linear case).
\end{remark}

\begin{remark}
In \cite{ReckeNonper} are proven results 
of the type of Theorems \ref{main1} and \ref{main2}
for 
periodic homogenization with localized defects of semilinear ODE systems of the type
$$
\Big(\big(A(\x)+B(\x)\big)u'(x)+c(x,u(x))\Big)'=d(x,u(x))\mbox{ for } x \in (0,1)
$$
with periodic $A \in L^\infty(\R;\M_n)$ and localized defects $B \in L^\infty(\R;\M_n)\cap L^1(\R;\M_n)$.
And in \cite{I} are proven results 
of the type of Theorems \ref{main1} and \ref{main2}
for
periodic homogenization (without localized defects) of quasilinear ODE systems of the type
$$
a(x,\x,u(x),u'(x))'=b(x,\x,u(x),u'(x))
\mbox{ for } x \in (0,1),
$$
where the maps $a(x,\cdot,u,u')$ and $a(x,\cdot,u,u')$ are supposed to be continuous from $\R$ into $\R^n$ and
periodic. 
\end{remark}

\begin{remark}
\label{Neuk}
What concerns existence and local uniqueness
for nonlinear periodic homogenization problems
(without assumption of global uniqueness), besides \cite{I,II,ReckeNonper}
we know only the result 
\cite{Bun}  for scalar semilinear elliptic PDEs of the type
$
\nabla\cdot(a(\x) 
\nabla u(x))=f(x)g(u(x)),
$
where the nonlinearity $g$ is supposed to have a sufficiently small local Lipschitz constant (on an appropriate bounded interval). Let us mention also \cite{Lanza1,Lanza2,Riva}, where existence and local uniqueness for  periodic homogenization problems for the linear Poisson equation with highly oscillating nonlinear Robin boundary conditions is shown. There the specific structure of the problem (no highly oscillating diffusion coefficients) allows to apply the classical implicit function theorem.
\end{remark}

\begin{remark}
\label{cub}
Below we will use several times
(in \reff{square}, \reff{e1ii} and \reff{rve})
the following simple fact: Denote by $\|\cdot\|$ the Euclidean norm in $\R^N$. Then
\bee
\label{cu}
\sup_{r>0,\,\ve>0,\, x \in \R^N}(r+\ve)^{-N}\int_{\|\xi-x\|<r} w(\xi/\ve)d\xi<\infty
\mbox{ for all $\Z^N$-periodic }w \in L^1_{\rm loc}(\R^N).
\ee
Indeed, for $r>0$, $\ve>0$, $x\in \R^N$ 
and $w \in L^1_{\rm loc}(\R^N)$ we have
$$
\int_{\|\xi-x\|<r} w(\xi/\ve)d\xi
=\ve^N\int_{\|y\|<r/\ve} w(x/\ve+y)dy
\le \ve^N\int_{(-r/\ve,r/\ve)^N} w(x/\ve+y)dy.
$$
But the cube $(-r/\ve,r/\ve)^N$ is covered by the cube $(-[r/\ve]-1,[r/\ve]+1)^N$
of integer edge length $2([r/\ve]+1)$ (where $[r/\ve]$ is the integer part of $r/\ve$),
and this cube can be covered by the union of $(2([r/\ve]+1))^N$ cubes of edge length one. Therefore the $\Z^N$-periodicity  of the function $w$ yields
$$
\ve^N\int_{(-r/\ve,r/\ve)^N} w(x/\ve+y)dy\le 
\ve^N(2([r/\ve]+1))^N \int_{(0,1)^N}w(y)dy\le \mbox{const}\, (r+\ve)^N,
$$
where the constant does not depend on $r$,
$\ve$ and $x$.
\end{remark}

~\\ 

Our paper is organized as follows: 
In Section \ref{secabstract} we consider abstract nonlinear parameter depending equations of the type
\bee
\label{intrabstract}
F_\ve(u)=0.
\ee
Here $\ve>0$ is the parameter. We prove 
Theorem \ref{ift}, which is
a result on existence and local uniqueness of a family of solutions $u=u_\ve \approx \bar u_\ve$ to \reff{intrabstract} with $\ve \approx 0$, where $\bar u_\ve$ is a family of  approximate solutions to \reff{intrabstract}, i.e. a family with
$F_\ve(\bar u_\ve)\to 0$ for $\ve \to 0$, and we estimate the norm of the error $u_\ve-\bar u_\ve$ by the norm of the discrepancy $F_\ve(\bar u_\ve)$. 
In the past this type of generalized implicit function theorems has been  applied to nonlinear singularly perturbed ODEs and elliptic and parabolic PDEs (see \cite{Butetc,But2022,Fiedler,
Magnus1,Magnus2,NURS,OR2009,OmelchenkoRecke2015,Recke2022,
ReckeOmelchenko2008}) as well as to  homogenization of nonlinear ODEs and PDEs
(see \cite{I,II, ReckeNonper}).
The proofs of all these results are based on the generalized implicit function theorem
of R. J. Magnus \cite[Theorem 1.2]{Magnus1} and on several of its modifications (see, e.g.
Theorem \ref{ift} below). Hence, these 
generalized implicit function theorems permit a general approach to nonlinear singular perturbation and  homogenization.

Contrary to the classical implicit function theorem, in Theorem \ref{ift} it is not supposed that
the linearized operators $F'_\ve(u)$ converge for $\ve \to 0$ in the uniform operator norm. And, indeed, in the applications to singularly perturbed problems or to  homogenization problems they do not converge for $\ve \to 0$ in the uniform operator norm. 
Remark that in the classical implicit function theorem one cannot omit, in general, the assumption, that $F_\ve'(u)$ converges for $\ve \to 0$ with respect to the uniform operator norm (cf. \cite[Section 3.6]{Katz}).

In place of a convergence for $\ve \to 0$ of the linearized operators $F'_\ve(u)$, in Theorem \ref{ift} it is assumed the weaker assumption \reff{Fas}, which allows to adapt the proof of the classical implicit function theorem (application of Banach's fixed point theorem) to this case. But, again contrary to the classical implicit function theorem, in the proof of Theorem \ref{ift} Banach's fixed point theorem has to be applied on closed balls with $\ve$-depending center points, and all estimates (radius of the balls, contraction constant of the mapping in the fixed point equation) must be shown to be $\ve$-independent.

In Sections \ref{sec3}--\ref{sec:proof2ii} we prove Theorems \ref{main1} and \ref{main2}, respectively, by means of 
the results of Section~\ref{secabstract}.
Here the main work is to construct 
appropriate families $\uu$ of approximate  solutions 
to \reff{BVP1} and \reff{BVP2}, respectively, 
with small 
(in appropriate function space norms)
discrepancies for $\ve \to 0$.
The reason, why in Theorem \ref{main1}(ii)  
better estimates for $\|u_\ve-u_0\|_\infty$ are possible,
under an  additional regularity assumption on $u_0$ ($u_0 \in W^{2,p_0}$ with certain $p_0$ larger than the space dimension), is easy to explain: By means of those more regular $u_0$ one can construct better families $\uu$, i.e. those with smaller discrepancies.

In order to apply implicit function theorems mainly one needs isomorphism properties 
of the linearized operators. In the setting of Sections \ref{sec3} and \ref{sec: proof1ii} for \reff{BVP1} they follow from K. Gr\"oger's results \cite{G}
about   maximal regularity of boundary value problems for elliptic equations and systems
with non-smooth data in  pairs of Sobolev spaces $W_0^{1,p}(\Omega;\R^n)$ and $W^{-1,p}(\Omega;\R^n)$ with $p\approx 2$
(cf. Theorem \ref{maxreg1} below).
In the setting of Sections \ref{sec4} and \ref{sec:proof2ii} for \reff{BVP2} they follow from the results
\cite[Lemma 6.2 and Theorem 6.3]{GR}
about   maximal regularity of boundary value problems for elliptic equations
with non-smooth data in  pairs of Sobolev-Morrey spaces $W_0^{1,p,\om}(\Omega)$ and $W^{-1,2,\om}(\Omega)^n$ with $\om \approx N-2$
(cf. Theorem \ref{maxreg2} below).

In the Theorems \ref{maxreg1} and 
Theorem \ref{maxreg2} below, which concern maximal
$W^{1,p}$ and $W^{1,2,\om}$ regularity of solutions to linear elliptic boundary value problems with non-smooth data, the allowed exponents $p$ and $\om$ depend on the non-smothness of the data. In general $p$ has to be close to (but larger than) two, and  $\om$ has to be close to (but larger than) $N-2$, respectively. 
But  this small increase of solution regularity is sufficient for proving that the localized defects $b^{\al \beta}_{ij}$ and $b_{ij}$ in \reff{BVP1} and \reff{BVP2}, respectively, are negligible (see Lemmas \ref{Blemma1} and \ref{Blemma2} below).

In order to apply implicit function theorems one needs also $C^1$-smoothness of the appearing nonlinear superposition operators. In the
settings of Sections
\ref{sec3}--\ref{sec:proof2ii}
these operators are well-defined and $C^1$-smooth on the used Sobolev spaces $W^{1,p}(\Omega;\R^n)$ with $p>2$ 
and on the used Sobolev-Morrey spaces 
$W^{1,2,\om}(\Omega)$
with $\om>N-2$, respectively, 
because these function spaces are continuously embedded into H\"older spaces. Here we use the approaches of \cite{GREvol,GrR,PRS,Recke1995}, where the classical implicit function theorem is applied to quasilinear elliptic and parabolic PDEs with non-smooth data (but without highly oscillating coefficients) by using maximal elliptic and parabolic regularity in pairs of Sobolev and Sobolev-Morrey spaces, respectively.

\section{An abstract result of implicit function theorem type}
\label{secabstract}
\setcounter{equation}{0}
\setcounter{theorem}{0}
Let $U$ and $V$ be  Banach spaces with norms $\|\cdot\|_U$ and $\|\cdot\|_V$, respectively.
For $\ve>0$ let be given 
$$
\uu \in U \mbox{ and } F_\ve\in C^1(U;V).
$$
We consider the abstract equation
\bee
\label{abeq}
F_\ve(u)=0. 
\ee
Roughly speaking, we will show the following:
If the elements $\uu$ approximately satisfy \reff{abeq} for $\ve \approx 0$, i.e. if $\|F_\ve(\uu)\|_V\to 0$ for  $\ve \to 0$, 
and if they are non-degenerate solutions
(cf. assumption \reff{coerz} below), then
for all $\ve \approx 0$ there exists exactly one solution $u=u_\ve$ to \reff{abeq} with $\|u-\uu\|_U\approx 0$, and $\|u_\ve-\uu\|_U=O(\|F_\ve(\uu)\|_V)$ for $\ve \to 0$.
For that we do not suppose 
that $\uu$ or $F_\ve(u)$ (with fixed $u \in U$) converge in $U$ for $\ve \to 0$ or that
$F'_\ve(u)$ (with fixed $u \in U$) converges in $\LL(U)$ for $\ve \to 0$
(cf. also Remark \ref{unb} below)
or that $\|\uu\|_U$ is bounded for $\ve \approx 0$.

\begin{theorem}
\label{ift}
Suppose that 
\bee
\label{coerz}
\left.
\begin{array}{l}
\mbox{there exist $\ve_0>0$ and $\rho>0$ such that for all $\ve \in (0,\ve_0]$ the operators $F_\ve'(\uu)$ are}\\
\mbox{Fredholm of index zero from }
$U$ \mbox{ into } $V$,
\mbox{and } \|F_\ve'(\uu))u\|_V \ge \rho \|u\|_U \mbox{ for all } u \in U
\end{array}
\right\}
\ee
and
\bee
\label{Fas}
\sup_{\|v\|_U \le 1}\|(F'_\ve(\uu+u)-F'_\ve(\uu))v\|_V \to 0
\mbox{ for } \ve + \|u\|_U \to 0
\ee
and
\bee
\label{lim}
\|F_\ve(\uu)\|_V \to 0 \mbox{ for } \ve \to 0.
\ee

Then there exist $\ve_1 \in (0,\ve_0]$ and $\delta>0$  such that for all $\ve \in (0,\ve_1]$ 
there exists exactly one  solution $u = u_\ve$ to~\reff{abeq} with $\|u-\uu\|_U \le \delta$. Moreover, for all $\ve \in (0,\ve_1]$ we have
\bee
\label{apriori1}
\|u_\ve-\uu\|_U \le \frac{2}{\rho}\|F_\ve(\uu))\|_V.
\ee
\end{theorem}
{\bf Proof }
Take $\ve \in (0,\ve_0]$. 
Because of assumption \reff{coerz} the linear operator $F'_\ve(\uu)$ is an isomorphism from $U$ onto $V$, and, hence,
equation \reff{abeq} is equivalent to the fixed point problem
$$
u=G_\ve(u):=u-F_\ve'(\uu)^{-1}F_\ve(u).
$$
Take $u_1,u_2 \in U$. Then 
\begin{eqnarray}
\label{strict}
\|G_{\ve}(u_1) - G_{\ve}(u_2)\|_U
&=&\left\|F_\ve'(\uu)^{-1}\int_0^1\left(
F'_\ve(\uu)-
F'_{\ve}(s u_1 + (1-s) u_2)\right)ds (u_1 - u_2)\right\|_U\nonumber\\
&\le& \frac{1}{\rho}\, \max_{0 \le s \le 1}\left\|\big(F'_\ve(\uu)-F'_{\ve}(s u_1 + (1-s) u_2)\big) (u_1 - u_2)\right\|_V.
\end{eqnarray}  
Here we used that \reff{coerz} yields that
$\rho\|F_\ve'(\uu)^{-1}v\|_U \le \|v\|_V$
for all $v \in V$.

Denote ${\cal B}^r_\ve:=\{u \in U:\; \|u-\uu\|_U \le r\}$.
If $u_1,u_2 \in {\cal B}^r_\ve$,
then also $su_1+(1-s)u_2 \in {\cal B}^r_\ve$
for all $s \in [0,1]$. Therefore  it follows from \reff{Fas} and \reff{strict} that there exist $r_0>0$ and $\ve_1 \in (0,\ve_0]$ such that for all $\ve \in (0,\ve_1]$
the maps $G_\ve$ are strictly contractive with contraction constant $1/2$ on the closed balls ${\cal B}^{r_0}_\ve$.
Moreover, for all $\ve \in (0,\ve_1]$ and $u \in {\cal B}^{r_0}_\ve$ we have
\begin{eqnarray*}
\label{in}
\left\|G_{\ve}(u) - \uu\right\|_U &\le & \left\|G_{\ve}(u) - G_\ve(\uu)\right\|_U
+\left\|G_{\ve}(\uu) - \uu\right\|_U\\
&\le & \frac{r_0}{2}+\left\|F_\ve'(\uu))^{-1}F_\ve(\uu)\right\|_U
\le \frac{r_0}{2}+\frac{1}{\rho}\left\|F_\ve(\uu)\right\|_V,
\end{eqnarray*}
and \reff{lim} yields that $G_\ve$ maps ${\cal B}^{r_0}_\ve$ into itself if $\ve_1$ is taken sufficiently small.

Now, Banach's fixed point principle yields the existence and uniqueness assertions of Theorem \ref{ift}, and the estimate \reff{apriori1}
follows as above:
$$
\|u_\ve-\uu\|_U \le \|G_{\ve}(u_\ve) - G_\ve(\uu)\|_U+\|G_{\ve}(\uu) - \uu\|_U
\le \frac{1}{2}\|u_\ve-\uu\|_U
+\frac{1}{\rho}\left\|F_\ve(\uu)\right\|_V.
$$
\qed

\begin{remark}
In \cite{Butetc,But2022,Fiedler,
Magnus1,Magnus2,NURS,
OmelchenkoRecke2015,Recke2022,
ReckeOmelchenko2008}) various  versions of Theorem \ref{ift} are presented
and applied to nonlinear singularly perturbed problems. 
These versions differ slightly according to which problems they are applied (ODEs or elliptic or parabolic PDEs, stationary or time-periodic solutions, semilinear or quasilinear problems, smooth or nonsmooth data, $\ve$-depending norms in $U$ and $V$).

For other results of the type of Theorem \ref{ift} and their applications to semilinear elliptic PDEs with numerically  determined approximate solutions  see \cite[Theorem 2.1]{Breden}
and \cite{Ca}.
\end{remark}

If one applies Theorem \ref{ift}, for example to boundary value problems for elliptic PDEs, then different choices of function spaces $U$ and $V$ and of their norms $\|\cdot\|_U$ and $\|\cdot\|_V$ 
and of the family $\uu$ of approximate solutions
are appropriate. Criteria for these choices often are the following: The family $\uu$  should be "simple" (for example, $\uu$ should be $\ve$-independent or given more less explicit in closed formulas, or to determine $ \uu$ numerically should be much cheeper than 
to determine the exact solution $u_\ve$ numerically), and the rate of convergence to zero of $\|F_\ve(\uu)\|_V$ for $\ve \to 0$ should be high. The norm $\|\cdot\|_V$ should be weak and the norm $\|\cdot\|_U$ should be strong such that the error estimate 
\reff{apriori1}
is strong. But at the same time 
the norm $\|\cdot\|_U$ should be weak such that the domain of local uniqueness, which contains all $u \in U$ with $\|u-\uu\|_U\le \delta$, is large.
These criteria are contradicting, of course. Hence, in any application of Theorem \ref{ift}
the choices of $U$, $V$, $\|\cdot\|_U$, $\|\cdot\|_V$ and $\uu$ are compromises according to the requirements of the application.

For example, in the present paper we will work with different choices of the spaces $U$ and  $V$ and their norms $\|\cdot\|_U$ and $\|\cdot\|_V$ and of
the families $\uu$:
Spaces and norms \reff{UVdef1} with exponents \reff{pdef} and families \reff{barudef1}
for proving Theorem \ref{main1}(i),
\reff{UVdef1} with \reff{pdefii} and \reff{barudef1ii}
for proving Theorem~\ref{main1}(ii)
and 
\reff{UVdef2} with \reff{omdef} and \reff{barudef2}
for proving Theorem \ref{main2}.
Remark that all the families $\uu$, we will work with, do not depend on the localized defects $b_{ij}^{\al \beta}$ and
$b_{ij}$.

One rather general way to find such compromises is described in  Corollary \ref{cor} below. It delivers existence and local uniqueness of solutions $u=u_\ve$ to the equation $F_\ve(u)=0$ with $\ve \approx 0$ and $\|u-u_0\|_\infty\approx 0$,
where $u_0 \in U$ is a given element
and $\|\cdot\|_\infty$ is another norm in $U$, which is allowed to be much weaker than the norm $\|\cdot\|_U$. 
The price for that is that the estimate \reff{apriori2} below of the error $u_\ve-u_0$ is with respect to the weaker norm $\|\cdot\|
_\infty$, only.
In Corollary \ref{cor} we use the notation $\|\cdot\|_\infty$ for a second  norm in  $U$, because in most of the applications to ODEs or PDEs (as well as in our applications, see Sections 
\ref{sec3}--\ref{sec:proof2ii}) it is an  $L^\infty$-norm. 

Remark that the space $U$ is supposed to be complete with respect to the norm $\|\cdot\|_U$, but not with respect to the norm $\|\cdot\|_\infty$, in general, and that it is not supposed that $\|F_\ve(u_0)\|_V$  
converges to zero for $\ve \to 0$.

\begin{corollary}
\label{cor}
Let be given $u_0 \in U$ such that
\bee
\label{coerza}
\left.
\begin{array}{l}
\mbox{there exist $\ve_0>0$ and $\rho>0$ such that for all $\ve \in (0,\ve_0]$ the operators $F_\ve'(u_0)$ are}\\
\mbox{Fredholm of index zero from }
$U$ \mbox{ into } $V$,
\mbox{and } \|F_\ve'(u_0))u\|_V \ge \rho \|u\|_U \mbox{ for all } u \in U,
\end{array}
\right\}
\ee
and let be given a norm $\|\cdot\|_\infty$ in $U$ such that
\begin{eqnarray}
\label{weaker}
&& \sigma:=\sup\{\|u\|_\infty: \; u \in U,\, \|u\|_U \le 1\}< \infty,\\
\label{Fas1}
&&\sup_{\|v\|_U \le 1}\|(F'_\ve(u_0+u)-F'_\ve(u_0))v\|_V \to 0
\mbox{ for } \ve + \|u\|_\infty \to 0.
\end{eqnarray}
Further, suppose that
\bee
\label{newconv}
\|\uu-u_0\|_\infty+\|F_\ve(\uu)\|_V \to 0 \mbox{ for } \ve \to 0.
\ee

Then
there exist $\ve_1 \in (0,\ve_0]$ and $\delta>0$  such that for all $\ve \in (0,\ve_1]$ 
there exists exactly one  solution $u = u_\ve$ to~\reff{abeq} with $\|u-u_0\|_\infty \le \delta$,
and 
\bee
\label{apriori2}
\|u_\ve-u_0\|_\infty \le \|\uu-u_0\|_\infty+
\frac{4\sigma}{\rho}\|F_\ve(\uu))\|_V.
\ee
\end{corollary}
{\bf Proof }
For all $\ve>0$ and $u,v \in U$ we have
\begin{eqnarray*}
&&\|(F_\ve'(\uu+u)-F_\ve'(\uu))v\|_V\\
&&\le
\|(F_\ve'(u_0+(\uu-u_0+u))-F_\ve'(u_0))v\|_V+
\|(F_\ve'(u_0)-F_\ve'(u_0+(\uu-u_0)))v\|_V.
\end{eqnarray*}
Hence, assumptions 
\reff{Fas1} and 
\reff{newconv} imply that
\bee
\label{newest}
\sup_{\|v\|_U \le 1}\|(F_\ve'(\uu+u)-F_\ve'(\uu))v\|_V
\mbox{ for } \ve+\|u\|_\infty \to 0.
\ee
In particular, \reff{Fas} is satisfied
(because of 
\reff{weaker}).
Similarly, for all $\ve>0$ and $u \in U$ we have that
$\|(F_\ve'(\uu)u\|_V
\ge
\|(F_\ve'(u_0)u\|_V-
\|(F_\ve'(u_0)-F_\ve'(u_0+(\uu-u_0)))u\|_V$.
Hence, \reff{coerza} and \reff{Fas1} yield that there exists $\ve_1 \in (0,\ve_0]$ such that 
$$
\|(F_\ve'(\uu)u\|_V \ge \frac{\rho}{2}\|u\|_U \mbox{ for all } \ve \in (0,\ve_1]
\mbox{ and } u \in U.
$$
Therefore \reff{coerz} is satisfied (with 
$\rho$ replaced by $\rho/2$).
Hence, assumptions \reff{coerza},
\reff{Fas1} and 
\reff{newconv} imply 
that \reff{coerz} is satisfied
(with $\ve_1$ in place of $\ve_0$ and with $\rho/2$ in place of $\rho$ in \reff{coerz}).
Hence, Theorem \ref{ift} yields the existence assertion of Corollary \ref{cor} and the error estimate
$$
\|u_\ve-u_0\|_\infty \le \|\uu-u_0\|_\infty
+\sigma\|u_\ve-\uu\|_U\le \|\uu-u_0\|_\infty+\frac{4\sigma}{\rho}\|F_\ve(\uu)\|_V.
$$

Now let us prove the local uniqueness assertion of Corollary \ref{cor}.
Take $\ve \in (0,\ve_1]$ and a solution $u \in U$ to \reff{abeq}. Then
$$
0=F_\ve(u)=F_\ve(\uu)+F'_\ve(\uu)(u-\uu)+
\int_0^1\left(F'_\ve(su+(1-s)\uu)-F'_\ve(\uu)\right)(u-\uu)ds,
$$
i.e.
$$
u-\uu=-F'_\ve(\uu)^{-1}\left(
F_\ve(\uu)+\int_0^1\left(F'_\ve(su+(1-s)\uu)-F'_\ve(\uu)\right)(u-\uu)ds\right),
$$
i.e.
\bee
\label{unest}
\|u-\uu\|_U\le \frac{1}{\rho}\left(
\|F_\ve(\uu)\|_V+\max_{0\le s \le 1}\|\left(F'_\ve(su+(1-s)\uu)-F'_\ve(\uu)\right)(u-\uu)\|_V\right).
\ee
But \reff{newest} yields that
$$
\max_{0\le s \le 1}\|\left(F'_\ve(su+(1-s)\uu)-F'_\ve(\uu)\right)(u-\uu)\|_V
=o(\|u-\uu\|_U) \mbox{ for }
\ve +\|u-\uu\|_\infty \to 0.
$$
Therefore \reff{newconv} implies that
$$
\max_{0\le s \le 1}\|\left(F'_\ve(su+(1-s)\uu)-F'_\ve(\uu)\right)(u-\uu)\|_V
=o(\|u-\uu\|_U) \mbox{ for }
\ve +\|u-u_0\|_\infty \to 0.
$$
Hence, if $\ve$ and $\|u-u_0\|_\infty$ are sufficiently small, then \reff{unest} yields that $\|u-\uu\|_U$ is small also, and the local uniqueness assertion of Theorem \ref{ift} implies that $u=u_\ve$.
\qed

\begin{remark}
\label{unb}
Assumption \reff{newconv} of Corollary \ref{cor} implies that
$\|\uu\|_\infty=O(1)$ for $\ve \to 0$,
and, hence, conclusion \reff{apriori2} yields that also $\|u_\ve\|_\infty=O(1)$ for $\ve \to 0$.
But in many applications of Theorem \ref{ift} and Corollary \ref{cor} one 
has $\|u_\ve\|_U\to \infty$ for $\ve \to 0$, and, therefore, one has
to use
families $\uu$ of approximate solutions with 
$\|\uu\|_U\to \infty$ for $\ve \to 0$
(cf. \reff{apriori1})!
For example, this is the case, in general,
in Sections \ref{sec3} and \ref{sec4}, where Corollary \ref{cor} is applied to the boundary value problems \reff{BVP1}
and \reff{BVP2}, respectively.
\end{remark}

\begin{remark}
In most of the applications of Corollary \ref{cor} to PDEs the element $u_0$ and the norm $\|\cdot\|_\infty$ are a priori given, and one has to choose  Banach spaces $U$ and $V$ 
and their norms $\|\cdot\|_U$ and
$\|\cdot\|_V$ such that 
\reff{weaker} is satisfied and
that the PDE problem is equivalent to an abstract equation $F_\ve(u)=0$ with 
$F_\ve\in C^1(U;V)$ such that 
\reff{coerza} and \reff{Fas1} are satisfied,
and one has to construct (by means of ansatzes, stretched variables, formal expansions in powers of $\ve$, cut-off functions, smoothing operators etc.) a family $\bar u_\ve$ with \reff{newconv}.
But at the beginning one does not know if existence and local uniqueness for $\ve \approx 0$ and $\|u-u_0\|_\infty\approx 0$ is true or not for the given PDE problem. If not, then one is trying to choose and to construct something, which does not exist.

For example, in Theorems \ref{main1} and \ref{main2}, which are the results of  applications of Corollary \ref{cor} to 
the boundary value problems \reff{BVP1} and \reff{BVP2}, respectively, the spaces $U$ and $V$ and their norms and the families $\bar u_\ve$ are hidden only, they do not appear in the formulations of Theorems \ref{main1} and \ref{main2}.
\end{remark}

\section{Proof of Theorem  \ref{main1}(i)}
\label{sec3}
\setcounter{equation}{0}
\setcounter{theorem}{0}
In this section we will prove Theorem \ref{main1}(i)
by means of Corollary \ref{cor}. 
For that we use the objects of Theorem \ref{main1}: The bounded Lipschitz domain $\Omega \subset \R^2$, the diffusion coefficients $a^{\al \beta}_{ij}$ and $b^{\al \beta}_{ij}$ with \reff{aass1}-\reff{monb1}, the drift and reaction functions $c_i^{\al},d^\al:\Omega\times \R^n\to \R$ with
\reff{diffass1}, 
the periodic correctors $v^{\beta}_j \in W^{1,2}_{\rm loc}(\R^2;\R^n)$, which are defined by the cell problems \reff{cell1},
the homogenized diffusion coefficients
$\hat a^{\al \beta}_{ij}\in \R$, which are
defined in \reff{hatAdef1}, and
 the
weak solution $u_0\in W_0^{1,2}(\Omega;\R^n)
\cap L^\infty(\Omega;\R^n)$ to the homogenized boundary value problem
\reff{hombvp1}.

\subsection{Maximal Sobolev regularity for linear elliptic systems}
\label{Sobolev}

As usual, the norm in the Sobolev space $W^{1,p}(\Omega;\R^n)$ (with $p\ge 2$) is denoted by 
$$
\|u\|_{1,p}:=
\left(\sum_{\al=1}^n
\int_\Omega
\left(|u^\al(x)|^p+
\sum_{i=1}^2 
|\partial_{x_i}u^\al(x)|^p\right)
dx\right)^{1/p}.
$$
The subspace $W_0^{1,p}(\Omega;\R^n)$ is the closure with respect to this norm 
of the set 
of all $C^{\infty}$-maps $u:\Omega \to \R^n$ with compact support,
and $W^{-1,p}(\Omega;\R^n):=W_0^{1,p'}(\Omega;\R^n)^*$ is the dual space to 
$W_0^{1,p'}(\Omega;\R^n)$ with 
$1/p+1/p'=1$
and with dual space norm
$$
\|\phi\|_{-1,p}:=\sup\{\langle \phi,\vp\rangle_{1,p'}:\; \vp \in W_0^{1,p'}(\Omega;\R^n),
\|\vp\|_{1,p'}\le1\},
$$
where
$
\langle \cdot,\cdot\rangle_{1,p'}:
W^{-1,p}(\Omega;\R^n)\times
W_0^{1,p'}(\Omega;\R^n)\to \R
$
is the dual pairing. 

Further, we introduce linear bounded operators $\hat A:W^{1,2}(\Omega;\R^n)\to W^{-1,2}(\Omega;\R^n)$ and, for $\ve>0$,
$A_\ve,B_\ve:W^{1,2}(\Omega;\R^n)\to W^{-1,2}(\Omega;\R^n)$
by
\bee
\label{ABdef1}
\left.
\begin{array}{l}
\displaystyle\langle \hat A u,\vp\rangle_{1,2}:=
\int_\Omega \hat a_{ij}^{\al \beta}
\partial_{x_j}u^\beta(x)\partial_{x_i}\vp^\al(x)dx,\\
\displaystyle\langle A_\ve u,\vp\rangle_{1,2}:=
\int_\Omega a_{ij}^{\al \beta}(\x)
\partial_{x_j}u^\beta(x)\partial_{x_i}\vp^\al(x)dx,\\
\displaystyle\langle B_\ve u,\vp\rangle_{1,2}:=
\int_\Omega b_{ij}^{\al \beta}(\x)
\partial_{x_j}u^\beta(x)\partial_{x_i}\vp^\al(x)dx,
\end{array}
\right\}
\mbox{ for all } \vp \in W_0^{1,2}(\Omega;\R^n).
\ee
Because of assumption \reff{aass1} and the H\"older inequality we have the following:
For any $p \ge 2$ the restrictions of $\hat A$, $A_\ve$ and $B_\ve$ to
$W^{1,p}(\Omega;\R^n)$
map $W^{1,p}(\Omega;\R^n)$ into
$W^{-1,p}(\Omega;\R^n)$, and
$$
\|\hat A u\|_{-1,p}+
\|A_\ve u\|_{-1,p}
+\|B_\ve u\|_{-1,p}
\le \mbox{const}\;\|u\|_{1,p} \mbox{ for all } u \in W^{1,p}(\Omega;\R^n),
$$
where the constant does not depend on $\ve$, $p$ and $u$. 
Further, because of assumption \reff{mona1} 
we have that 
$$
\inf\left\{\langle A_\ve u,u\rangle_{1,2}:\;
\ve>0,\;
u \in W_0^{1,2}(\Omega;\R^n),
\; \|u\|_{1,2}=1\right\}>0,
$$
and similarly for $A_\ve +B_\ve$ and for $\hat A$.
Therefore K. Gr\"oger's maximal regularity results for elliptic systems with non-smooth data \cite[Theorems 1 and 2 and Remark 14]{G} imply the following:
\begin{theorem}
\label{maxreg1}
There exist $p_1>2$ and $\rho>0$ such that
for all $\ve>0$ and all $p\in[2,p_1]$
the linear operators $\hat A$, $A_\ve$ and $A_\ve+B_\ve$
are bijective from $W_0^{1,p}(\Omega;\R^n)$ onto $W^{-1,p}(\Omega;\R^n)$ and that
\bee
\label{Ainvert1}
\|\hat A^{-1}\phi\|_{1,p}+
\|A_\ve^{-1}\phi\|_{1,p}+ 
\|(A_\ve+B_\ve)^{-1}\phi\|_{1,p} \le \rho\, \|\phi\|_{-1,p}
\mbox{ for all } \phi \in W^{-1,p}(\Omega;\R^n).
\ee
\end{theorem}

\begin{remark}
Estimates of the type \reff{Ainvert1} often are called Meyers' estimates
because of the initiating paper \cite{M}
of N.G. Meyers, see also \cite{Ga,ME1975}. For the case of smooth boundaries $\partial \Omega$ see 
\cite[Theorem 4.1]{Ben}. See also \cite{CP} for $\Omega$ being a cube and for continuous diffusion coefficients as well as for certain transmission problems and applications to linear elliptic periodic homogenization.
\end{remark}
\begin{remark}
\label{remsing}
It is easy to verify that the linear operators $A_\ve$ do not converge for $\ve \to 0$ with respect to the uniform operator norm in ${\cal L}(W^{1,p}(\Omega;\R^n);W^{-1,p}(\Omega;\R^n))$ for certain $p\ge 2$, in general (see \cite[Remark 8.4]{Ci} and Lemma
\ref{Klemma} below). This is the reason why the classical implicit function theorem is not directly applicable to the boundary value problem \reff{BVP1}, in general. But the linear operators $A_\ve$ converge in a certain weak sense to
the linear operator
$\hat A$ for $\ve \to 0$ (cf. Theorem \ref{Shentheorem1} and Lemma \ref{Klemma}
below), 
and this is the reason why Theorem \ref{ift} is applicable to \reff{BVP1}.
\end{remark}

\begin{remark}
Even if the diffusion coefficients $a_{ij}^{\al \beta}$ would be constant, i.e. if 
$a_{ij}^{\al \beta}(y)=\hat a_{ij}^{\al \beta}$ for all $y \in \R^2$, i.e. if $A_\ve=\hat A$ for all $\ve$, then the classical implicit function theorem would not be applicable to the boundary value problem \reff{BVP1}, because the linear operators $B_\ve$ are not small in $\LL(W^{1,p}(\Omega);W^{-1,p}(\Omega))$ for small $\ve$, in general. They are small only in $\LL(W^{1,p_1}(\Omega);W^{-1,p}(\Omega))$ with $p_1>p$, cf. Lemma \ref{Blemma1} below.
\end{remark}

The question of $W^{1,p}$-regularity of the correctors $v^{\beta}_j$ is much simpler than Theorem \ref{maxreg1} because there do not appear problems with boundary regularity.
It is well-known (cf., e.g. \cite[Chapter 2.2, formula (2.2.22)]{Shen}) that 
\bee 
\label{vsmooth1}
\mbox{there exists } p_2>2 \mbox{ such that }
v^{\beta}_j\in W_{\rm loc}^{1,p_2}(\R^2;\R^n), \mbox{ in particular }
v^{\beta}_j\in L^{\infty}(\R^2;\R^n).
\ee

\subsection{Abstract setting for the boundary value problem \reff{BVP1}}
\label{abstract1}
Now we introduce an abstract setting of the type  of Corollary \ref{cor} for the boundary value problem~\reff{BVP1}.  We take 
$p_1$ from  Theorem \ref{maxreg1}  
and $p_2$ from \reff{vsmooth1} and fix
exponents $p$ and $p'$ as follows:
\bee
\label{pdef}
2<p< \min\{p_1,p_2\},
\; p':=\frac{p}{p-1}.
\ee
The Banach spaces $U$ and $V$  and their norms are defined by
\bee
\label{UVdef1}
U:=W_0^{1,p}(\Omega;\R^n),\;
V:=W^{-1,p}(\Omega;\R^n),\;
\|\cdot\|_U:=\|\cdot\|_{1,p},\;
\|\cdot\|_V:=\|\cdot\|_{-1,p}.
\ee
The second norm $\|\cdot\|_\infty$ in $U$ of Corollary \ref{cor} is defined by \reff{normdef}.
Because the space dimension is two, the assumption \reff{weaker} of Corollary \ref{cor} is satisfied in this setting.
Further,  the $C^1$-smooth operators $F_\ve:U \to V$ of Theorem \ref{ift} are defined by 
$$
F_\ve(u):=(A_\ve+B_\ve) u+C(u)
$$
with a 
nonlinear operator $C:
L^\infty(\Omega;\R^n)\to W^{-1,2}(\Omega;\R^n)$ defined by
\bee
\label{Cdef1}
\langle C(u),\vp\rangle_{1,2}:=
\int_\Omega\Big(c_i^\al(x,u(x))\partial_{x_i}\vp^\al(x)+d^\al(x,u(x))\vp^\al(x)\Big)dx
\mbox{ for all } \vp \in W_0^{1,2}(\Omega;\R^n).
\ee
Because of assumption \reff{diffass1} we have that  
the nonlinear operator $C$ is $C^1$-smooth from $L^\infty(\Omega;\R^n)$ into $W^{-1,2}(\Omega;\R^n)$,
and
$$
\langle C'(u)v,\vp\rangle_{1,2}:=
\int_\Omega\Big(\partial_{u^\gamma}c_i^\al(x,u(x))\partial_{x_i}\vp^\al(x)+\partial_{u^\gamma}d^\al(x,u(x))\vp^\al(x)\Big)v^{\gamma}(x)dx.
$$
Remark that
the linear operator $C$ can be considered as a compact operator from $W^{1,p}(\Omega;\R^n)$ into $W^{-1,p}(\Omega;\R^n)$
because the space $W^{1,p}(\Omega;\R^n)$
is compactlyly embedded into the space
$L^\infty(\Omega;\R^n)$ (because the dimension of $\Omega$ is two).
In particular, the linear operators
$A_\ve+B_\ve+C'(u)$ are Fredholm of index zero from $W_0^{1,p}(\Omega;\R^n)$ into $W^{-1,p}(\Omega;\R^n)$ because of Theorem~\ref{maxreg1}.

With these choices a  vector function $u$ is a weak solution to the boundary value problem \reff{BVP1} if and only if $u$ belongs to the function space $U$ and satisfies the operator equation
$F_\ve(u)=0$.
Here we used Theorem \ref{maxreg1} again.

Finally, we define the  element $u_0\in U$ of Corollary \ref{cor} to be the solution $u_0$ of the homogenized boundary value problem \reff{hombvp1}, which is given in Theorem \ref{main1}.
Remark that 
\bee
\label{unulleq1}
\hat Au_0+C(u_0)=0.
\ee
Therefore Theorem \ref{maxreg1} yields that 
\bee
\label{usmooth1}
u_0 \in W_0^{1,p_1}(\Omega;\R^n),
\mbox{ in particular }
u_0 \in L^{\infty}(\Omega;\R^n),
\ee
i.e. $u_0 \in U=W_0^{1,p}(\Omega;\R^n)$, as needed.

In order to prove Theorem \ref{main1}(i) we have to verify the conditions \reff{coerza} and \reff{Fas1}
of  Corollary \ref{cor}
in the setting introduced above, i.e. that
there exists $\ve_0>0$ such that
\bee
\label{coerz1}
\inf\left\{\|(A_\ve+B_\ve +C'(u_0))u\|_{-1,p}:\; \ve \in (0,\ve_0], 
 u \in W_0^{1,p}(\Omega;\R^n),
 \|u\|_{1,p}=1\right\}>0
\ee
and that
\bee
\label{Cas1}
\sup_{\|v\|_{1,p} \le 1}
\|(C'(u_0+u)-C'(u_0))v\|_{-1,p}
\to 0 \mbox{ for } 
\|u\|_\infty \to 0,
\ee
and we have to construct a family $\bar u_\ve \in W_0^{1,p}(\Omega;\R^n)$ such that  \reff{newconv} is
satisfied in the setting introduced above, i.e. that 
\bee
\label{lim1}
\|\uu-u_0\|_\infty +\|(A_\ve+B_\ve)\uu+C(\uu)\|_{-1,p}
\to 0 \mbox{ for } 
\ve \to 0.
\ee
This is what we are going to do below in the next three subsections.

\subsection{Verification of  \reff{coerz1}}
\label{subcoerz1}
Suppose that \reff{coerz1} is not true. Then there exist sequences $\ve_1,\ve_2,\ldots>0$ and $u_1,u_2,\ldots \in W_0^{1,p}(\Omega;\R^n)$ such that 
\bee
\label{conv1}
\lim_{l\to \infty}\Big(\ve_l+\|(A_{\ve_l}+B_{\ve_l}+C'(u_0))u_l\|_{-1,p}\Big)=0,
\ee
but
\bee
\label{norm1}
\|u_l\|_{1,p}= 1 \mbox{ for all } l.
\ee
Because $W^{1,p}(\Omega;\R^n)$ is reflexive and because it is compactly embedded into  $L^\infty(\Omega;\R^n)$, without loss of generality we may assume that there exists $u_*\in W^{1,p}(\Omega;\R^n)$ such that 
\bee
\label{infconv1}
u_l\rightharpoonup u_* \mbox{ for $l\to \infty$
weakly in }
W^{1,p}(\Omega;\R^n)
\mbox{ and }
\lim_{l \to \infty}\|u_l-u_*\|_\infty=0.
\ee
The functions $u_l$ are continuous and vanish on the boundary $\partial \Omega$,
and for $l\to \infty$ they converge uniformly on $\overline \Omega$ to the continous function $u_*$.
Therefore $u_*\in W_0^{1,p}(\Omega;\R^n)$.

Moreover, $C'(u_0)$ is a linear bounded operator from $L^\infty(\Omega;\R^n)$ into $W^{-1,p}(\Omega;\R^n)$, therefore \reff{conv1} and \reff{infconv1} imply that
\bee
\label{conver1}
\lim_{l\to \infty}\|(A_{\ve_l}+B_{\ve_l})u_l+C'(u_0)u_*\|_{-1,p}=0.
\ee

The next lemma is the place where we use that the coefficients $b_{ij}^{\al \beta}$ are localized defects in the sense of assumption \reff{aass1}, i.e. that they belong to $L^1(\R^2)$:
\begin{lemma}
\label{Blemma1}
For any $r>s\ge 2$ we have
$$
\sup\left\{\|B_{\ve}u\|_{-1,s}:\;
u \in W^{1,r}(\Omega;\R^n),\; \|u\|_{1,r} 
\le 1\right\}=O\left(\ve^{2/s-2/r}\right) \mbox{ for } \ve \to 0.
$$
\end{lemma}
{\bf Proof }
Because of $b^{\al \beta}_{ij} \in L^\infty(\R^2)\cap
L^1(\R^2)$ (cf. assumption \reff{aass1}) we have for any $t \in [1,\infty)$ that
\bee
\label{best}
\int_{\R^2}\left|b^{\al \beta}_{ij}(y)\right|^tdy=
\left\|b^{\al \beta}_{ij}\right\|^t_{L^\infty(\R^2)}\int_{\R^2}
\left(\frac{\left|b^{\al \beta}_{ij}(y)\right|}{\left\|b^{\al \beta}_{ij}\right\|_{L^\infty(\R^2)}}\right)^tdy
\le \left\|b^{\al \beta}_{ij}\right\|^{t-1}_{L^\infty(\R^2)}
\left\|b^{\al \beta}_{ij}\right\|_{L^1(\R^2)}.
\ee
Now, take $r>s\ge 2$ and define $s'\in [2,\infty)$ and $t\in (2,\infty)$ by
$1/s+1/s'=1$ and $1/s-1/r=1/t$.
Then $1/r+1/s'+1/t=1$,
and the generalized H\"older inequality implies for $u \in W^{1,r}(\Omega)$ that
\begin{eqnarray*}
\|B_{\ve}u\|_{-1,s}&=&
\sup_{\|\vp\|_{1,s'}\le 1}
\int_\Omega b^{\al \beta}_{ij}(\x)\partial_{x_j}u^\beta(x)\partial_{x_i}\vp^\al(x)dx\\
&\le& \|u\|_{1,r}\left(\sum_{i,j=1}^2\sum_{\al,\beta=1}^n\int_\Omega |b^{\al \beta}_{ij}(\x)|^tdx\right)^{1/t}\\
&\le&\ve^{2/t}\;
\|u\|_{1,r}\left(\sum_{i,j=1}^2\sum_{\al,\beta=1}^n
\int_{\R^2} |b^{\al \beta}_{ij}(y)|^tdy\right)^{1/t}.
\end{eqnarray*}
Hence, 
\reff{best}
implies the assertion of the lemma.
\qed\\

If we apply Lemma \ref{Blemma1} with $r=p$ and $s=2$ (cf. \reff{pdef}), then  \reff{norm1} implies that $\|B_{\ve_l}u_l\|_{-1,2}\to 0$ for $l \to \infty$. Hence, \reff{conver1} yields that
\bee
\label{convert1}
\lim_{l\to \infty}\|A_{\ve_l}u_l
+C'(u_0)u_*\|_{-1,2}=0.
\ee

Now we use the following theorem, which is well-known in periodic homogenization theory for linear elliptic equations and  systems with $L^\infty$-coefficients (see, e.g. \cite[Lemma 8.6]{Che}, \cite[Theorem 2.3.2]{Shen}, \cite{Tartar}). 
Its proof is based on the Div-Curl Lemma.
It claims that the sequence of diffusion tensors $a_{ij}^{\al \beta}(\cdot/\ve_1),a_{ij}^{\al \beta}(\cdot/\ve_2),\ldots$
H-converges to the diffusion tensor $\hat a_{ij}^{\al \beta}$. We formulate this in the language of the differential operators
$A_{\ve_l}$ and $\hat A$.
\begin{theorem}
\label{Shentheorem1}
Let be given $u \in W^{1,2}(\Omega;\R^n)$ and $\phi\in W^{-1,2}(\Omega;\R^n)$ and sequences $\ve_1,\ve_2,\ldots>0$ and $u_1,u_2,\ldots \in 
W^{1,2}(\Omega;\R^n)$ such that
$$
u_l\rightharpoonup u \mbox{ for $l\to \infty$ weakly in 
$W^{1,2}(\Omega;\R^n)$ and }
\lim_{l \to \infty}\left(\ve_l+\|A_{\ve_l}u_l-\phi\|_{-1,2}\right)=0.
$$
Then $\hat A u=\phi$.
\end{theorem}

Because of \reff{infconv1}, \reff{convert1} and Theorem \ref{Shentheorem1} it follows that
$
(\hat A+C'(u_0))u_*=0,
$ 
i.e. that
$u_*$ is a weak solution
to the linearized boundary value problem
\reff{linhombvp1}. Hence, by assumption of Theorem~\ref{main1}, we get that $u_*=0$.
Therefore
\reff{conver1} implies that $\|(A_{\ve_l}+B_{\ve_l})u_l\|_{-1,p}\to 0$ for $l \to \infty$.
But this contradicts to \reff{Ainvert1} and \reff{norm1}.

\subsection{Verification of  \reff{Cas1}}
\label{sub3a}
Because of assumption \reff{diffass1} the functions $u \mapsto \partial_{u_\gamma}c^\al_i(\cdot,u)$
and $u \mapsto \partial_{u_\gamma}d^\al(\cdot,u)$
are continuous from $\R^n$ into $L^\infty(\Omega)$ and, hence, uniformly continuous on bounded sets. Therefore
$$
[\phi^{\al \gamma}_i(u)](x):=
\partial_{u_\gamma}c^\al_i(x,u_0(x)+u(x))-\partial_{u_\gamma}c^\al_i(x,u_0(x))
$$
and 
$$
[\psi^{\al \gamma}(u)](x):=
\partial_{u_\gamma}d^\al(x,u_0(x)+u(x))-\partial_{u\gamma}d^\al(x,u_0(x))
$$
tend to zero for $\|u\|_\infty \to 0$
uniformly with respect to $x \in \Omega$.
Moreover, for $u,v \in L^\infty(\Omega;\R^n)$
and $\vp\in W_0^{1,p'}(\Omega;\R^n)$
we have
\begin{eqnarray*}
&&\langle(C'(u_0+u)-C'(u_0))v,\vp\rangle_{1,p'}=\int_\Omega\Big(
[\phi_i^{\al \gamma}(u)](x)
\partial_{x_i}\vp^\al(x)
+[\psi^{\al \gamma}(u)](x)
\vp^\al(x)\Big)v^\gamma(x) dx\\
&&\le \mbox{const}\sum_{\al,\gamma=1}^n\left(
\sum_{i=1}^2
\|\phi^{\al \gamma}_i(u)\|_\infty
+\|\psi^{\al \gamma}(u)\|_\infty\right)
\|v\|_{\infty}\|\vp\|_{1,p'},
\end{eqnarray*}
where the constant does not depend on $u$, $v$ and $\vp$.
Hence, \reff{Cas1} is proved.

\begin{remark}
In fact, we proved slightly more than \reff{Cas1}, namly that $\langle(C'(u_0+u)-C'(u_0))v,\vp\rangle_{1,p'}$ tends to zero for $\|u\|_\infty \to 0$ not only uniformly with respect to $\|v\|_{1,p}\le 1$ and $\|\vp\|_{1,p'}\le 1$, but even 
uniformly with respect to $\|v\|_{\infty}\le 1$ and $\|\vp\|_{1,p'}\le 1$.
\end{remark}

\subsection{Construction of approximate solutions with \reff{lim1}}
\label{sublim1}
In this subsection we will construct a family $\uu\in  W_0^{1,p}(\Omega;\R^n)$ 
with 
\reff{lim1}.
For that we will do some calulations which are well-known in periodic homogenization theory for linear elliptic problems (cf., e.g. \cite[Chapter 3.2]{Shen}), and some estimates which seem to be new. 

We denote by $\|\cdot\|$ the Euclidean norm in $\R^2$ and take a mollifier function, i.e. a $C^\infty$-function $\rho:\R^2\to \R$ such that
$$
\rho(x)\ge 0 \mbox{ and } \rho(x)=\rho(-x) \mbox{ for all } x \in \R^2,\;
\rho(x)=0 \mbox{ for } \|x\| \ge 1,
\int_{\R^2}\rho(x)dx=1,
$$
and for $\delta>0$ we define linear Steklov  smoothing operators $S_\delta: L^1(\Omega) \to C^\infty(\R^2)$ by
$$
[S_\delta u](x):=\int_{\Omega}\rho_\delta(x-\xi) u(\xi)d\xi
\mbox{ with }
\rho_\delta(x):=\rho(x/\delta)/\delta^2.
$$

\begin{lemma}
\label{Slemma}
For all $r > 1$ we have
\begin{eqnarray}
\label{Sconv}
&&\lim_{\delta\to 0}\int_\Omega|u(x)-[S_\delta u](x)|^rdx=0
\mbox{ for all } u \in L^r(\Omega),\\
\label{pest1a}
&&\sup\left\{\delta^{2}\|S_\delta u\|^r_\infty:\;\delta \in (0,1],\;  u \in L^{r}(\Omega),\;
\int_\Omega|u(x)|^rdx\le 1 
\right\}< \infty,\\
\label{pest2a}
&&\sup\left\{\delta^r\int_\Omega|\partial_{x_i}S_\delta u(x)|^rdx:\;\delta \in (0,1],\;  u \in L^{r}(\Omega),\;
\int_\Omega|u(x)|^rdx\le 1, 
\right\}<\infty.
\end{eqnarray}
\end{lemma}
{\bf Proof }
Assertion \reff{Sconv} is proved e.g. in
\cite[Lemma 1.1.1]{Barbu}.

In order to prove assertions \reff{pest1a} and \reff{pest2a}  take $\delta>0$, $r,r'>1$ with $1/r+1/r'=1$, and take $u \in L^r(\Omega)$.
Then the H\"older inequality implies that for all $x\in \Omega$ we have 
$$
|[S_\delta u](x)|=\left|\int_{\Omega} u(\xi)\rho_\delta(x-\xi)^{1/r}\rho_\delta(x-\xi)^{1/r'}d\xi\right|\le\left(\int_{\Omega}|u(\xi)|^r\rho_\delta(x-\xi)d\xi\right)^{1/r}.
$$
Here we used that
$
\int_{\R^2}\rho_\delta(x-\xi)d\xi
=
\int_{\R^2}\rho(\xi)d\xi=1.
$
It follows that
$$
|[S_\delta u](x)|^r
\le \frac{1}{\delta^2}
\int_{\Omega}|u(\xi)|^r\rho((x-\xi)/\delta)d\xi
\le \mbox{const } \frac{1}{\delta^2}\int_{\Omega}|u(\xi)|^rd\xi,
$$
where the constant does not depend on $\delta$, $x$ and $u$.
Hence, \reff{pest1a} is proved.

In order to prove assertion \reff{pest2a} we use that
$\partial_{x_i}\rho_\delta(x)=-\partial_{x_i}\rho(x/\delta)/\delta^3$ for all $x \in \Omega$ and, therefore, that
$$
\int_{\R^2}|\partial_{x_i}\rho_\delta(x-\xi)|d\xi=\frac{1}{\delta^3}\int_{\R^2}|\partial_{x_i}\rho((x-\xi)/\delta)|d\xi=
\frac{1}{\delta}\int_{\R^2}|\partial_{x_i}\rho(\xi)|d\xi.
$$
It follows that
\begin{eqnarray*}
|[\partial_{x_i}S_\delta u](x)|&\le&
\left(\int_{\Omega}| u(\xi)|^r|\partial_{x_i}\rho_\delta(x-\xi)|d\xi\right)^{1/r}
\left(\int_{\Omega}|\partial_{x_i}\rho_\delta(x-\xi)|d\xi\right)^{1/r'}\\
&\le& \mbox{ const }
\frac{1}{\delta^{1/r'}}\left(\int_{\Omega}|u(\xi)|^r|\partial_{x_i}\rho_\delta(x-\xi)|d\xi\right)^{1/r}
\end{eqnarray*}
and, hence,
\begin{eqnarray*}
\int_\Omega\left|[\partial_{x_i}S_\delta u](x)\right|^rdx&\le& \mbox{ const }
\frac{1}{\delta^{r/r'}}\int_{\Omega}| u(\xi)|^r\int_{\Omega}|\partial_{x_i}\rho_\delta(x-\xi)|dx\;d\xi\\
&\le&\mbox{ const }
\frac{1}{\delta^r}\int_{\Omega}|u(\xi)|^rd\xi,
\end{eqnarray*}
where, again, the constants do not depend on $\delta$, $x$ and $u$.
Hence, \reff{pest2a} is proved also.
\qed\\

Further, for $\ve>0$ we define the boundary strip 
$\Omega_\ve:=\left\{x\in \Omega:\; 
\inf_{y \in \partial\Omega}
\|x-y\|<\ve\right\}$ of size $\ve$.
It follows that
\bee
\label{Omest}
|\Omega_\ve|=O(\ve)
\mbox{ for } \ve \to 0,
\ee
where $|\Omega_\ve|$ is the two-dimensional Lebesque measure
of $\Omega_\ve$.
Further, we take a family $\eta_\ve$ of
cut-of functions of size $\ve$, i.e. of $C^\infty$-functions $\Omega\to \R$ such that
\bee
\label{etaest}
\left.
\begin{array}{l}
0 \le \eta_\ve(x)\le 1 \mbox{ for } x \in \Omega,\\
\eta_\ve(x)=1\mbox{ for } x \in  \Omega\setminus \Omega_{2\ve},\\
\eta_\ve(x)=0\mbox{ for } x \in \Omega_{\ve},\\
\sup\left\{\ve\,|\partial_{x_i}\eta_\ve(x)|:\; \ve>0, \;x \in \Omega,\; i=1,2\right\}<\infty.
\end{array}
\right\}
\ee

Now we define the needed family $\uu \in W_0^{1,p}(\Omega;\R^n)$, which should satisfy
\reff{lim1}. For $\ve>0$, $x \in \Omega$ and $\al=1,\ldots,n$ we set
\bee
\label{barudef1}
\uu^\al(x):=
u^\al_0(x)+\ve \eta_\ve(x) [S_{\delta_\ve}\partial_{x_k}u_0^\gamma](x) v_k^{\gamma \al}(\x)
\mbox{ with } \delta_\ve:=\frac{1}{|\ln \ve|}.
\ee
Remark that
\reff{vsmooth1}
and \reff{usmooth1} yield that
\bee
\label{p*def}
\uu \in W_0^{1,p_*}(\Omega;\R^n)
\mbox{ with } p_*:=\min\{p_1,p_2\},
\ee
and \reff{pdef} implies that $\uu\in W_0^{1,p}(\Omega;\R^n)$, as needed.
Further, from \reff{vsmooth1}, \reff{usmooth1}
and \reff{pest1a} and from definition 
\reff{barudef1} it follows that
\bee
\label{Kest1}
\|\uu-u_0\|_\infty \le \mbox{const }\frac{\ve}{\delta_\ve^{2/p_1}}\|u_0\|_{1,p_1}\le \mbox{const }\ve |\ln \ve|^{2/p_1},
\ee
where the constants do not depend on $\ve$.
Hence,
\reff{lim1} is satisfied if 
$\|(A_\ve+B_\ve)\uu+C(\uu)\|_{-1,p}\to 0$ for $\ve \to 0$. But the operator $C$ is continuous from $L^\infty(\Omega;\R^n)$ into $W^{-1,p}(\Omega;\R^n)$, therefore \reff{Kest1} yields that it remains to show that $\|(A_\ve+B_\ve)\uu+C(u_0)\|_{-1,p}\to 0$ for $\ve \to 0$. Finally, using \reff{unulleq1} we get that it remains to show that 
$$
\|(A_\ve+B_\ve)\uu-\hat Au_0\|_{-1,p}\to 0 \mbox{ for } \ve \to 0.
$$
But this is the assertion of
the two lemmas below.
\begin{lemma}
\label{BLemma1}
We have 
$\|B_\ve\uu\|_{-1,p}\to 0$ for $\ve \to 0$.
\end{lemma}
{\bf Proof }
Because of Lemma \ref{Blemma1}  with $r=p_*$ and $s=p$ 
(cf. \reff{pdef} and \reff{p*def}) we have that
\bee
\label{Buuest}
\|B_\ve\uu\|_{-1,p}=O(\ve^{2/p-2/p_*}\|\uu\|_{1,p_*}) \mbox{ for } \ve \to 0.
\ee
Hence, 
it sufficies to show that 
$\ve^{2/p-2/p_*}\|\uu\|_{1,p_*}\to 0$
for $\ve \to 0$ or, what is the same,
that
$$
\lim_{\ve\to 0}
\ve^{2/p-2/p_*}\left(\int_\Omega |\partial_{x_j}\big(\uu^\al(x)-u_0^\al(x)\big)|^{p_*}dx\right)^{1/p_*}=0.
$$
But $\partial_{x_j}\big(\uu^\al(x)-u_0^\al(x)\big)$ equals to
$$
\eta_\ve(x) [S_{\delta_\ve}\partial_{x_k}u^\gamma_0](x) \partial_{y_j}v^{\gamma \al}_k(\x)
+\ve\Big(
\partial_{x_j}\eta_\ve(x) 
[S_{\delta_\ve}\partial_{x_k}u_0](x) + \eta_\ve(x) 
[S_{\delta_\ve}\partial_{x_j}\partial_{x_k}u_0](x)\Big) v_k(\x)
$$
(cf. \reff{barudef1}), therefore we get the following: 
It remains to prove that 
\begin{eqnarray}
\label{e1}
&&\lim_{\ve \to 0}\ve^{2/p-2/p_*}
\left(\int_\Omega |\eta_\ve(x) [S_{\delta_\ve}\partial_{x_k}u^\gamma_0](x) \partial_{y_j}v^{\gamma \al}_k(\x)|^{p_*}dx\right)^{1/p_*}=0,\\
\label{e2}
&&\lim_{\ve \to 0}\ve^{1+2/p-2/p_*}
\left(\int_\Omega |\partial_{x_j}\eta_\ve(x) [S_{\delta_\ve}\partial_{x_k}u^\gamma_0](x)v^{\gamma \al}_k(\x)|^{p_*}dx\right)^{1/p_*}=0,\\
\label{e3}
&&\lim_{\ve \to 0}\ve^{1+2/p-2/p_*}
\left(\int_\Omega |\eta_\ve(x) [\partial_{x_j}S_{\delta_\ve}\partial_{x_k}u^\gamma_0](x)v^{\gamma \al}_k(\x)|^{p_*}dx\right)^{1/p_*}=0.
\end{eqnarray}

First, let us prove \reff{e1}. Because of  the regularity \reff{usmooth1} of $u_0$
and because of \reff{pest2a} with
$r=p_*=\min\{p_1,p_2\}$
 we have
$$
|\eta_\ve(x) [S_{\delta_\ve}\partial_{x_k}u^\gamma_0](x) \partial_{y_j}v^{\gamma \al}_k(\x)|
\le \mbox{const } \frac{1}{\delta_\ve^{2/p_*}}
\sum_{\gamma=1}^n\sum_{k=1}^2
|\partial_{y_j}v^{\gamma \al}_k(\x)|,
$$
where the constant does not depend on $x$ and $\ve$. 

Take $r>0$ such that $\Omega \subseteq [-r,r]^2$. Because of 
\reff{vsmooth1} the $\Z^2$-periodic functions $\partial_{y_j}v^{\gamma \al}_k$ belong to $L_{loc}^{p_*}(\R^2)$. Therefore \reff{cu} yields that
\bee
\label{square}
\sup_{\ve \in (0,1]}\int_\Omega |\partial_{y_j}v^{\gamma \al}_k(\x)|^{p_*}dx\le\sup_{\ve \in (0,1]} \int_{(-r,r)^2} |\partial_{y_j}v^{\gamma \al}_k(\x)|^{p_*}dx<\infty.
\ee
Hence,
$$
\ve^{2/p-2/p_*}
\left(\int_\Omega |\eta_\ve(x) [S_{\delta_\ve}\partial_{x_k}u^\gamma_0](x) \partial_{y_j}v^{\gamma \al}_k(\x)|^{p_*}dx\right)^{1/p_*}\le 
\mbox{const }\ve^{2/p-2/p_*}\,|\ln \ve|^{2/p_*},
$$
where the constant does not depend on $\ve$.

Now, let us prove \reff{e2}. Because of
\reff{vsmooth1}  and \reff{etaest}
we have
$$
|\partial_{x_j}\eta_\ve(x) [S_{\delta_\ve}\partial_{x_k}u^\gamma_0](x) v^{\gamma \al}_k(\x)|
\le \mbox{const } \frac{1}{\ve}\,
|[S_{\delta_\ve}\partial_{x_k}u^\gamma_0](x)|,
$$
where the constant does not depend on $x$ and $\ve$. Hence, using \reff{usmooth1}
we get
$$
\ve^{1+2/p-2/p_*}
\left(\int_\Omega |\partial_{x_j}\eta_\ve(x) [S_{\delta_\ve}\partial_{x_k}u^\gamma_0](x) v^{\gamma \al}_k(\x)|^{p_*}dx\right)^{1/p_*}\le 
\mbox{const }\ve^{2/p-2/p_*},
$$
where the constant does not depend on $\ve$.

And finally, let us prove \reff{e3}.
Because of \reff{vsmooth1}, 
we have
$$
|\eta_\ve(x) [\partial_{x_j}S_{\delta_\ve}\partial_{x_k}u^\gamma_0](x) v^{\gamma \al}_k(\x)|
\le \mbox{const } 
\sum_{\gamma=1}^n\sum_{k=1}^2|\partial_{x_j}S_{\delta_\ve}\partial_{x_k}u^\gamma_0](x)|,
$$
where the constant does not depend on $x$ and $\ve$. Because of
\reff{usmooth1} and \reff{pest2a}
it follows that
$$
\ve^{1+2/p-2/p_*}
\left(\int_\Omega |\eta_\ve(x) [\partial_{x_j}S_{\delta_\ve}\partial_{x_k}u^\gamma_0](x) v^{\gamma \al}_k(\x)|^{p_*}dx\right)^{1/p_*}\le 
\mbox{const }\ve^{1+2/p-2/p_*}
|\ln\ve|,
$$
where the constant does not depend on $\ve$.
\qed\\

\begin{remark}
For a result related to Lemma \ref{BLemma1} see \cite[Theorem 7.4]{BDL}.
\end{remark}

\begin{lemma}
\label{Klemma}
We have
$\|A_\ve\uu-\hat Au_0\|_{-1,p}\to 0$
for $\ve \to 0$.
\end{lemma}
{\bf Proof }
We 
insert the definition \reff{barudef1} of $\uu$ into $\langle A_\ve\uu-\hat Au_0,\vp\rangle_{1,p'}$
with $\vp \in C^\infty(\Omega;\R^n)$
and calculate as follows:
\begin{eqnarray}
&&
\langle A_\ve\uu-\hat Au_0,\vp\rangle_{1,p'}
\nonumber\\
&&=\int_\Omega\left(
a^{\al \beta}_{ij}(\x)\partial_{x_j}\left(u_0^\beta+\ve\eta_\ve
[S_{\delta_\ve}\partial_{x_k}u_0^\gamma]
v_k^{\beta \gamma}(\x)
\right)                                    
-\hat a^{\al \beta}_{ij}\partial_{x_j}
u_0^\beta\right)
\partial_{x_i}\vp^\al dx\nonumber\\
&&=
\int_\Omega\left(
\left(a^{\al \beta}_{ij}(\x)
-\hat a^{\al \beta}_{ij}\right)
\partial_{x_j}u_0^\beta
+a^{\al \beta}_{ij}(\x)
\eta_\ve
[S_{\delta_\ve}\partial_{x_k}u_0^\gamma]
\partial_{y_j}v_k^{\beta \gamma}(\x)
\right)                                    
\partial_{x_i}\vp^\al dx\nonumber\\
&&\;\;\;\;\;\;+\ve\int_\Omega
a^{\al \beta}_{ij}(\x)\partial_{x_j}(
\eta_\ve(x)
[S_{\delta_\ve}\partial_{x_k}u_0^\gamma](x))
v_k^{\beta \gamma}(\x)
\partial_{x_i}\vp^\al(x)
dx\nonumber\\
&&=
\int_\Omega\left(
a^{\al \beta}_{ij}(\x)
-\hat a^{\al \beta}_{ij}
+a^{\al \gamma}_{ik}(\x)
\partial_{y_k}v_j^{\gamma\beta}(\x)\right)\eta_\ve(x)
[S_{\delta_\ve}\partial_{x_j}u_0^\beta](x)
\partial_{x_i}\vp^\al(x)dx\nonumber\\
&&\;\;\;\;\;\;+\int_\Omega
\left(a^{\al \beta}_{ij}(\x)
-\hat a^{\al \beta}_{ij}\right)\left(
\partial_{x_j}u_0^\beta(x)-
\eta_\ve(x)
[S_{\delta_\ve}\partial_{x_j}u_0^\beta](x)\right)
\partial_{x_i}\vp^\al(x)dx\nonumber\\
&&\;\;\;\;\;\;+\ve\int_\Omega
a^{\al \beta}_{ij}(\x)\partial_{x_j}(
\eta_\ve(x)
[S_{\delta_\ve}\partial_{x_k}u_0^\gamma](x))
v_k^{\beta \gamma}(\x)
\partial_{x_i}\vp^\al(x)
dx.
\label{int1}
\end{eqnarray}

Further, for $\al,\beta=1,\ldots,n$ and $i,j,k=1,2$
we define
$\Z^2$-periodic auxiliary functions $f_{ij}^{\al \beta}\in L_{\rm loc}^{p}(\R^2)$ and $g_{ij}^{\al \beta}\in W_{\rm loc}^{2,p}(\R^2)$
and 
$h_{ijk}^{\al \beta}\in W_{\rm loc}^{1,p}(\R^2)$  
by
\begin{eqnarray}
\label{bdef1}
&&f_{ij}^{\al \beta}(y):=a_{ij}^{\al \beta}(y)
+a^{\al \gamma}_{ik}(y)
\partial_{y_k}v_j^{\gamma \beta}(y)-\hat a^{\al \beta}_{ij},\\
\label{cdef1}
&&\Delta g_{ij}^{\al \beta}(y)=f_{ij}^{\al \beta}(y),\; \int_{(0,1)^2}g_{ij}^{\al \beta}(y)dy=0,\\
\label{phidefa}
&&h_{ijk}^{\al \beta}(y):=
\partial_{y_i}g_{jk}^{\al \beta}(y)
-\partial_{y_j}g_{ik}^{\al \beta}(y).
\end{eqnarray}
The functions $h_{ijk}^{\al \beta}$  sometimes are called flux correctors, and, because they are $\Z^2$-periodic and locally $W^{1,p}$-functions with $p>2$, we have that
\bee
\label{hbound}
h_{ijk}^{\al \beta} \in L^\infty(\R^2).
\ee 
From \reff{hatAdef1} and \reff{bdef1} follows that $\int_{[0,1]^2}f^{\al \beta}_{ij}(y)dy=0$,
therefore problem \reff{cdef1} is uniquely strongly solvable with respect to $g^{\al \beta}_{ij}$. Further,
from \reff{cell1} follows that $\partial_{y_i}f_{ij}^{\al \beta}=0$.
Hence, \reff{cdef1} implies that $\partial_{y_i}g_{ij}^{\al \beta}=0$.
Therefore \reff{cdef1} and \reff{phidefa} yield that
\bee
\label{phiprop1}
\partial_{y_i}h_{ijk}^{\al \beta}=f_{jk}^{\al \beta} 
\mbox{ and }
h_{ijk}^{\al \beta}=-h_{kji}^{\al \beta}.
\ee
Using \reff{phiprop1} we get
\bee
\label{ibp1}
\ve \partial_{x_k}\Big(h_{ijk}^{\al \beta}(\x)\partial_{x_i}\vp^\beta(x)\Big)=
f_{ij}^{\al \beta}(\x)\partial_{x_i}\vp^\beta(x)
\mbox{ for all } \vp \in C^\infty(\Omega;\R^n)
\ee 
(this is \cite[formula (3.1.5)]{Shen}).

Now we insert \reff{bdef1} and \reff{ibp1} into \reff{int1}, integrate by parts, use $h_{kij}^{\al \beta}(\x)\partial_{x_k}\partial_{x_i}\vp^\al(x)=0$
(cf. \reff{phiprop1}) again,
and this way we get
\begin{eqnarray}
&&
\langle A_\ve\uu-\hat Au_0,\vp\rangle_{1,p'}
\nonumber\\
&&=\int_\Omega
f_{ij}^{\al \beta}(\x)
\eta_\ve(x)
[S_{\delta_\ve}\partial_{x_j}u_0^\beta](x)
\partial_{x_i}\vp^\al(x)
dx\nonumber\\
&&\;\;\;\;\;\;+\int_\Omega
\left(a^{\al \beta}_{ij}(\x)
-\hat a^{\al \beta}_{ij}\right)\left(
\partial_{x_j}u_0^\beta(x)-
\eta_\ve(x)
[S_{\delta_\ve}\partial_{x_j}u_0^\beta](x)\right)
\partial_{x_i}\vp^\al(x)dx\nonumber\\
&&\;\;\;\;\;\;+\ve\int_\Omega
a^{\al \beta}_{ij}(\x)\partial_{x_j}\Big(
\eta_\ve(x)
[S_{\delta_\ve}\partial_{x_k}u_0^\gamma](x)\Big)
v_k^{\beta \gamma}(\x)
\partial_{x_i}\vp^\al(x)
dx\nonumber\\
&&=\ve\int_\Omega\left(-h_{ijk}^{\al \beta}(\x)+
a^{\al\gamma}_{ik}(\x) v_j^{\gamma\beta}(\x)
\right)
\partial_{x_k}\Big(\eta_\ve(x)
[S_{\delta_\ve}\partial_{x_j}u_0^\beta](x)\Big)
\partial_{x_i}\vp^\al(x)dx\nonumber\\
&&\;\;\;\;\;\;+\int_\Omega
\left(a^{\al \beta}_{ij}(\x)
-\hat a^{\al \beta}_{ij}\right)\left(
\partial_{x_j}u_0^\beta(x)-
\eta_\ve(x)
[S_{\delta_\ve}\partial_{x_j}u_0^\beta](x)\right)
\partial_{x_i}\vp^\al(x)dx.
\label{threeint1}
\end{eqnarray}

Let us estimate  the right-hand side of \reff{threeint1}. 
In what follows we will use 
the fact that the functions $v^{\gamma \beta}_k$ and $h^{\al \beta}_{ijk}$ are bounded
(cf. \reff{vsmooth1} and \reff{hbound}), and all constants will be independent on $\ve$ and $\vp$.

Because of the H\"older inequality
we have
\begin{eqnarray*}
&&\left|\ve\int_\Omega\left(-h_{kij}^{\al \beta}(\x)+
a^{\al \gamma}_{ik}(\x) v_j^{\gamma \beta}(\x)
\right)\partial_{x_k}\eta_\ve(x)
[S_{\delta_\ve}\partial_{x_j}u_0^\beta](x)
\partial_{x_i}\vp^\al(x)dx\right|\nonumber\\
&&\le\mbox{const }\ve\left(\sum_{i=1}^2
\int_{\Omega_\ve}
|\partial_{x_i}\eta_\ve(x)|^pdx\right)^{1/p}
\sum_{\beta=1}^n\sum_{j=1}^2\|S_{\delta_\ve}\partial_{x_j}
u_0^\beta\|_\infty
\|\vp\|_{1,p'}\nonumber\\
&&\le\mbox{const } \frac{\ve^{1/p'}} 
{\delta_\ve^{2/p}}\|u_0\|_{1,p}
\|\vp\|_{1,p'}
\le
\mbox{const }\ve^{1/p'}|\ln \ve|^{2/p}
\|\vp\|_{1,p'}
\end{eqnarray*}
(here we used 
\reff{pest1a}
\reff{Omest} and \reff{etaest}) and
\begin{eqnarray*}
&&\left|\ve\int_\Omega\left(-h_{kij}^{\al \beta}(\x)+
a^{\al \gamma}_{ik}(\x) v_j^{\gamma \beta}(\x)
\right)\eta_\ve(x)
[\partial_{x_k}
S_{\delta_\ve}\partial_{x_j}u_0^\beta](x)
\partial_{x_i}\vp^\al(x)dx\right|\nonumber\\
&&\le\mbox{const }
\frac{\ve}{\delta_\ve}
\|u_0\|_{1,p}
\|\vp\|_{1,p'}\le\mbox{const }
\ve\, |\ln \ve|
\|\vp\|_{1,p'}
\end{eqnarray*}
(here we used \reff{pest2a})
and
\begin{eqnarray*}
&&\left|
\int_\Omega
\left(a^{\al \beta}_{ij}(\x)
-\hat a^{\al \beta}_{ij}\right)\left(
1-
\eta_\ve(x)\right)
[S_{\delta_\ve}\partial_{x_j}u_0^\beta](x)
\partial_{x_i}\vp^\al(x)dx\right|\nonumber\\
&&\le\mbox{const}\left(
\int_{\Omega_\ve}
|1-\eta_\ve(x)|^pdx\right)^{1/p}
\sum_{\beta=1}^n\sum_{j=1}^2\|[S_{\delta_\ve}\partial_{x_j}
u_0^\beta]\|_\infty
\|\vp\|_{1,q}\nonumber\\
&&\le\mbox{const } \frac{\ve^{1/p}}
{\delta_\ve^{2/p}}\|u_0\|_{1,p}
\|\vp\|_{1,p'}
\le\mbox{const } \ve^{1/p}|\ln \ve|^{2/p}
\|\vp\|_{1,p'}
\end{eqnarray*}
(here we used \reff{pest1a} and \reff{Omest}).
Further, we have
\begin{eqnarray}
\label{Sdeltaest1}
&&\left|
\int_\Omega
\left(a^{\al \beta}_{ij}(\x)
-\hat a^{\al \beta}_{ij}\right)\left(
\partial_{x_j}u_0^\beta(x)-
[S_{\delta_\ve}\partial_{x_j}u_0^\beta](x)\right)
\partial_{x_i}\vp^\al(x)dx\right|\nonumber\\
&&\le \mbox{const}\left(\sum_{\beta=1}^n\sum_{j=1}^2\int_\Omega\left|\partial_{x_j}u_0^\beta(x)-
[S_{\delta_\ve}\partial_{x_j}u_0^\beta](x)\right|^{p}dx\right)^{1/p}
\|\vp\|_{1,p'}.
\end{eqnarray}
But the right-hand side of \reff{Sdeltaest1} is $o(1)$ for $\ve \to 0$ uniformly with respect to $\|\vp\|_{1,p'}\le 1$ (cf. \reff{Sconv}).
Hence, the lemma is proved.
\qed\\

\begin{remark}
\label{boundaryterms}
Remark that in the proof above we did not need that the test function $\vp$ vanishes on the boundary $\partial \Omega$. In particular,
no boundary integrals appeared after the integration by parts 
in \reff{threeint1}
because
of the presence of the cut-off functions $\eta_\ve$ (no matter if $\vp$ vanishes on $\partial \Omega$ or not).
Therefore \reff{threeint1} can be used for proving results of the type of Theorem \ref{main1} for problems with other than Dirichlet boundary conditions also.
\end{remark}

\section{Proof of Theorem  \ref{main1}(ii)}
\label{sec: proof1ii}
\setcounter{equation}{0}
\setcounter{theorem}{0}
In this section we will prove Theorem \ref{main1}(ii)
by means of Corollary \ref{cor}. 
For that we use 
the assumption 
\bee
\label{unullass}
u_0 \in W^{2,p_0}(\Omega;\R^n)
\mbox{ with certain } p_0>2, \mbox{ in particular }
\partial_{x_k}u_0 \in L^\infty(\Omega;\R^n),
\ee
from Theorem \ref{main1}(ii)
and define the Banach spaces $U:=W_0^{1,p}(\Omega;\R^n)$  and $V:=W^{-1,p}(\Omega;\R^n)$ as in Section \ref{sec3}, but now with $p$ defined not in \reff{pdef}, but in
\bee
\label{pdefii}
2<p< p_*:=\min\{p_0,p_1,p_2\}.
\ee
And the family $\uu$ of approximate solutions to \reff{BVP1} is not defined by \reff{barudef1}, but by
\bee
\label{barudef1ii}
\uu^\al(x):=
u^\al_0(x)+\ve \eta_\ve(x) 
\partial_{x_k}u_0^\gamma(x) v_k^{\gamma \al}(\x).
\ee
Because of \reff{vsmooth1} and \reff{unullass} we have 
$\uu \in U$ and 
$\|\uu-u_0\|_\infty =O(\ve)$  for 
$\ve \to 0$.
Hence, in order to prove the assertions of Theorem \ref{main1}(ii) we have to show that  
\bee
\label{1ii1}
\mbox{ for certain } \la>0 \mbox{ we have }\|(A_\ve+B_\ve)\uu-\hat Au_0\|_{-1,p}=O(\ve^\la)
\mbox{ for } \ve \to 0
\ee 
and that
\bee
\label{1ii2}
\mbox{ for all } \la\in (0,1/2) \mbox{ we have }\|A_\ve\uu-\hat Au_0\|_{-1,p}=O(\ve^\la)
\mbox{ for } \ve \to 0.
\ee 
Both assertions \reff{1ii1} and
\reff{1ii2}
follow from the Lemmas
\ref{BLemma1ii} and \ref{Klemmaii} below.
\begin{lemma}
\label{BLemma1ii}
We have 
$\|B_\ve\uu\|_{-1,p}=O(\ve^{2/p-2/p_*})$ for $\ve \to 0$.
\end{lemma}
{\bf Proof } We proceed as in the proof of Lemma \ref{BLemma1}.  
Because of \reff{Buuest}
it sufficies to show that $\sup_{\ve \in (0,1]}\|\uu\|_{1,p_*}<\infty$.
But $u_0 \in W_0^{1,p_*}(\Omega;\R^n)$, therefore it sufficies to show that
$$
\sup_{\ve\in (0,1]}
\left(\int_\Omega |\partial_{x_j}(\uu^\al(x)-u_0^\al(x))|^{p_*}dx\right)^{1/p_*}<\infty.
$$
Because of
\begin{eqnarray*}
&&\partial_{x_j}\big(\uu^\al(x)-u_0^\al(x)\big)\\
&&=
\eta_\ve(x) \partial_{x_k}u^\gamma_0(x) \partial_{y_j}v^{\gamma \al}_k(\x)
+\ve\Big(
\partial_{x_j}\eta_\ve(x) 
\partial_{x_k}u^\gamma_0(x) + \eta_\ve(x) 
\partial_{x_j}\partial_{x_k}u^\gamma_0(x)\Big) v^{\gamma \al}_k(\x)
\end{eqnarray*}
we get the following: 
It remains to prove that 
\begin{eqnarray}
\label{e1ii}
&&\sup_{\ve\in (0,1]}
\left(\int_\Omega |\eta_\ve(x) \partial_{x_k}u^\gamma_0(x) \partial_{y_j}v^{\gamma \al}_k(\x)|^{p_*}dx\right)^{1/p_*}<\infty,\\
\label{e2ii}
&&\sup_{\ve\in (0,1]} \ve
\left(\int_\Omega |\partial_{x_j}\eta_\ve(x) \partial_{x_k}u^\gamma_0
(x)v^{\gamma \al}_k(\x)|^{p_*}dx\right)^{1/p_*}<\infty,\\
\label{e3ii}
&&\sup_{\ve\in (0,1]} \ve
\left(\int_\Omega |\eta_\ve(x) \partial_{x_j}\partial_{x_k}u^\gamma_0
(x)v^{\gamma \al}_k(\x)|^{p_*}dx\right)^{1/p_*}<\infty.
\end{eqnarray}

First, let us prove \reff{e1ii}. Because of   \reff{unullass} we have that
$$
|\eta_\ve(x) \partial_{x_k}u^\gamma_0(x) \partial_{y_j}v^{\gamma \al}_k(\x)|
\le \mbox{const } 
\sum_{\gamma=1}^n\sum_{k=1}^2
|\partial_{y_j}v^{\gamma \al}_k(\x)|,
$$
where the constant does not depend on $x$ and $\ve$. 
Hence, \reff{square} implies \reff{e1ii}.

Assertion \reff{e2ii} follows from 
$\ve \|\partial_{x_j}\eta_\ve\|_\infty \le \mbox{const}$ (cf.
\reff{etaest}) and from $v^{\gamma \al}_k \in L^\infty(\R^2)$
(cf.
\reff{vsmooth1}),
and, similarly, \reff{e3ii} follows from
$$
\ve
\left(\int_\Omega |\eta_\ve(x) [\partial_{x_j}\partial_{x_k}u^\gamma_0
(x) v^{\gamma \al}_k(\x)|^{p_*}dx\right)^{1/p_*}\le 
\mbox{const }\ve \sum_{k=1}^2 \left(
\int_\Omega|\partial_{x_j}\partial_{x_k}u^\gamma_0(x)|^{p_*}dx\right)^{1/p_*},
$$
where the constants do not depend on $\ve$.
\qed

\begin{lemma}
\label{Klemmaii}
We have
$\|A_\ve\uu-\hat Au_0\|_{-1,p}
=O(\ve^{1/p})$
for $\ve \to 0$.
\end{lemma}
{\bf Proof }
Take $\ve>0$ and $\vp \in W_0^{1,p'}(\Omega;\R^n)$.
As in \reff{threeint1} it follows that
\begin{eqnarray}
&&
\langle A_\ve\uu-\hat Au_0,\vp\rangle_{1,p'}
\nonumber\\
&&=\ve\int_\Omega\left(-h_{ijk}^{\al \beta}(\x)+
a^{\al\gamma}_{ik}(\x) v_j^{\gamma\beta}(\x)
\right)
\partial_{x_k}\Big(\eta_\ve(x)
\partial_{x_j}u_0^\beta(x)\Big)
\partial_{x_i}\vp^\al(x)dx\nonumber\\
&&\;\;\;\;\;\;+\int_\Omega
\left(a^{\al \beta}_{ij}(\x)
-\hat a^{\al \beta}_{ij}\right)\left(
\partial_{x_j}u_0^\beta(x)-
\eta_\ve(x)
\partial_{x_j}u_0^\beta(x)\right)
\partial_{x_i}\vp^\al(x)dx.
\label{threeint1ii}
\end{eqnarray}
Let us estimate  the right-hand side of \reff{threeint1ii}. We have
\begin{eqnarray*}
&&\left|\ve\int_\Omega\left(-h_{kij}^{\al \beta}(\x)+
a^{\al \gamma}_{ik}(\x) v_j^{\gamma \beta}(\x)
\right)\partial_{x_k}\eta_\ve(x)
\partial_{x_j}u_0^\beta(x)
\partial_{x_i}\vp^\al(x)dx\right|\nonumber\\
&&\le\mbox{const }\ve\left(\sum_{i=1}^2
\int_{\Omega_\ve}
|\partial_{x_i}\eta_\ve(x)|^pdx\right)^{1/p}
\|u_0\|_{1,\infty}
\|\vp\|_{1,p'}\le
\mbox{const }\ve^{1/p}
\|\vp\|_{1,p'}
\end{eqnarray*}
(here we used \reff{Omest} and \reff{etaest}) and
\begin{eqnarray*}
&&\left|\ve\int_\Omega\left(-h_{kij}^{\al \beta}(\x)+
a^{\al \gamma}_{ik}(\x) v_j^{\gamma \beta}(\x)
\right)\eta_\ve(x)
\partial_{x_k}
\partial_{x_j}u_0^\beta(x)
\partial_{x_i}\vp^\al(x)dx\right|\nonumber\\
&&\le\mbox{const }
\ve
\|u_0\|_{2,p}
\|\vp\|_{1,p'}\le\mbox{const }
\ve
\|\vp\|_{1,p'}
\end{eqnarray*}
and
\begin{eqnarray*}
&&\left|
\int_\Omega
\left(a^{\al \beta}_{ij}(\x)
-\hat a^{\al \beta}_{ij}\right)\left(
1-
\eta_\ve(x)\right)
\partial_{x_j}u_0^\beta(x)
\partial_{x_i}\vp^\al(x)dx\right|\nonumber\\
&&\le\mbox{const}\left(
\int_{\Omega_\ve}
|1-\eta_\ve(x)|^pdx\right)^{1/p}
\|
u_0\|_{1,\infty}
\|\vp\|_{1,p'}
\le\mbox{const } \ve^{1/p}|
\|\vp\|_{1,p'}
\end{eqnarray*}
(here we used \reff{Omest} again),
where the constants do not depend on $\ve$
and $\vp$.
Hence, the lemma is proved.
\qed\\

Because of \reff{pdefii} the exponent $p>2$ may be choosen arbitrarily close to two. Therefore
$\la=\min\{1/p,1/p'\}>1/2$ can be made arbitrarily close to $1/2$.
Hence, \reff{1ii2} is proved.

\section{Proof of Theorem  \ref{main2}(i)}
\label{sec4}
\setcounter{equation}{0}
\setcounter{theorem}{0}
In this section we will prove Theorem \ref{main2}(i)
by means of Corollary \ref{cor}. 
We proceed as in Section~\ref{sec3}.
We use the objects of Theorem \ref{main2}: The bounded Lipschitz domain $\Omega \subset \R^N$, the diffusion coefficients $a_{ij}$ and $b_{ij}$ with \reff{aass2}-\reff{monb2}, the drift and reaction functions $c_i,d:\Omega\times \R\to \R$ with
\reff{diffass2}, 
the periodic correctors $v_i \in W^{1,2}_{\rm loc}(\R^N)$, which are defined by the cell problems \reff{cell2},
the homogenized diffusion coefficients
$\hat a_{ij} \in \R$, which are
defined in \reff{hatAdef2}, 
and the
weak solution $u_0$ to the homogenized boundary value problem~\reff{hombvp2}.
By assumption of Theorem \ref{main2} we have
\bee
\label{pom}
u_0 \in W^{2,p_0}(\Omega)\cap 
W_0^{1,2}(\Omega)
 \mbox{ with certain } p_0>N,
\mbox{ in particular }
u_0,\partial_{x_i}u_0 \in L^{\infty}(\Omega).
\ee

\subsection{Maximal Sobolev-Morrey regularity for linear elliptic equations}
\label{Campanato}
In order to introduce Morrey spaces and Sobolev-Morey spaces (cf., e.g. \cite{BenF,Chen,Gia,Kufner,Sa,Troi}) we denote for $x \in \Omega$ and $r>0$
\bee
\label{Omxrdef}
\Omega_{x,r}:=\{\xi\in \Omega:\; \|\xi-x\|<r\},
\mbox{ where $\|\cdot\|$ is the Euclidean norm in } \R^N.
\ee
Take $\om \in [0,N)$.
The Morrey space $L^{2,\om}(\Omega)$ is the subspace of the Lebesgue $L^2(\Omega)$ with norm
$$
\|u\|_{2,\om}:=\sup_{x\in \Omega, r\in (0,1]}r^{-\om/2}\left(\int_{\Omega_{x,r}}|u(\xi)|^2d\xi\right)^{1/2} 
$$
which is defined by
$
L^{2,\om}(\Omega):=\{u\in L^2(\Omega):
\|u\|_{2,\om}<\infty\}$. 
The Morrey space $L^{2,\om}(\Omega)$
is continuously embedded into the Lebesgue space
$L^2(\Omega)$, and the Lebesgue space
$L^p(\Omega)$ 
is continuously embedded into the Morrey space
$L^{2,\om}(\Omega)$ for 
$\om= N(1-2/p)$.
In particular, \reff{pom} yields that
\bee
\label{pom1}
\partial_{x_i}u_0 \in W^{1,2,\om_0}(\Omega) \mbox{ with }
\om_0:=N\left(1-\frac{2}{p_0}\right).
\ee
The Sobolev-Morrey space $W^{1,2,\om}(\Omega)$ is the subspace of 
the Sobolev space  $W^{1,2}(\Omega)$ 
with norm
$$
\|u\|_{1,2,\om}:=\left(\int_\Omega|u(x)|^2dx+\sum_{i=1}^N\|\partial_{x_i}u\|^2_{2,\om}\right)^{1/2},
$$
which is defined by
$W^{1,2,\om}(\Omega):=\{u\in W^{1,2}(\Omega):\;
\|u\|_{1,2,\om}<\infty\}$.
It is well-known that
\bee
\label{embed}
\left.
\begin{array}{l}
W^{1,2,\om}(\Omega)
\mbox{ is compactly embedded into }
C^ {0,\al}(\overline \Omega)\\
\displaystyle\mbox{for } N-2<\om< N \mbox{ and } 0 \le \al <1-\frac{N-\om}{2}.
\end{array}
\right\}
\ee
Further, we denote $W_0^{1,2,\om}(\Omega):=W^{1,2,\om}(\Omega)\cap W_0^{1,2}(\Omega)$.

And finally, the Sobolev-Morrey space $W^{-1,2,\om}(\Omega)$ (sometimes called 
Sobolev-Morrey space of functionals)
is the subspace of $W^{-1,2}(\Omega)=W_0^{1,2}(\Omega)^*$ 
with norm 
$$
\|\phi\|_{-1,2,\om}:=\sup\left\{
r^{-\om/2}\langle \phi,\vp\rangle_{1,2}:\;
x\in \Omega,r\in (0,1],\vp \in W_0^{1,2}(\Omega),\|\vp\|_{1,2}\le 1, \mbox{supp} \,\vp \subset \Omega_{x,r}\right\},
$$
which is defined by
$W^{-1,2,\om}(\Omega):=\{\phi\in W^{-1,2}(\Omega):\;
\|\phi\|_{-1,2,\om}<\infty\}$, and it is continuously embedded into $W^{-1,2}(\Omega)$.
Here $\|\cdot\|_{1,2}$ is the norm in the Sobolev space $W^{1,2}(\Omega)$, and 
$\langle\cdot,\cdot\rangle_{1,2}:
W_0^{1,2}(\Omega)\times W^{-1,2}(\Omega)\to \R$
is the dual mapping (as introduced in Subsection \ref{Sobolev} for vecor functions).

It is well-known (cf. e.g. \cite[Theorem 3.9]{Griep})
that 
the linear map $(f,g)\mapsto \phi_{f,g}$, which is defined by
\bee
\label{phidef1}
\langle \phi_{f,g},\vp\rangle_{1,2}:=
\int_\Omega \Big(f_i(x)\partial_{x_i}\vp(x)+g(x)\vp(x)\Big)dx
\mbox{ for all } 
\vp \in 
W^{1,2}(\Omega),
\ee
is well-defined and bounded from $L^{2,\om}(\Omega)^N
\times L^{\infty}(\Omega)$
into $W^{-1,2,\om}(\Omega)$.
Therefore the linear operators $\hat A, A_\ve$ and $B_\ve$,
which are defined in \reff{ABdef1}, are bounded from 
$W^{1,2,\om}(\Omega)$ into $W^{-1,2,\om}(\Omega)$.
But moreover, for sufficiently large $\om$ 
the operators $\hat A,A_\ve$ and $A_\ve+B_\ve$ are
are isomorphisms
from $W_0^{1,2,\om}(\Omega)$ onto $W^{-1,2,\om}(\Omega)$:
\begin{theorem}
\label{maxreg2}
There exist $\om_1\in (N-2,N)$ and $\rho>0$ such that for all  $\ve>0$ and
$\om\in[\om_1,N)$ the linear operators 
$\hat A$, $A_\ve$ and $A_\ve+B_\ve$ 
are bijective from $W_0^{1,2,\om}(\Omega)$ onto $W^{-1,2,\om}(\Omega)$, 
and 
$$
\|\hat A^{-1}\phi\|_{1,2,\om}+
\|A_\ve^{-1}\phi\|_{1,2,\om} +
\|(A_\ve+B_\ve)^{-1}\phi\|_{1,2,\om}
\le \rho \|\phi\|_{-1,2,\om}
\mbox{ for all } \phi \in W^{-1,2,\om}(\Omega).
$$
\end{theorem}
Theorem \ref{maxreg2} follows from 
\cite[Theorem 4.1]{Griep} and
\cite[Lemma 6.2 and Theorem 6.3]{GR}, and its proof essentially uses 
the assumptions \reff{Omass2}-\reff{monb2}.

The question of Sobolev-Morrey regularity of the correctors $v_i$ is much simpler
than Theorem \ref{maxreg2} because there do not appear problems with boundary regularity. For example, the interior Sobolev-Morrey regularity result
\cite[Theorem 2.16]{Troi} yields that
\bee
\label{vreg2}
\mbox{there exists } \om_2\in(N-2,N)
\mbox{ such that }
v_i \in W_{\rm loc}^{1,2,\om_2}(\R^N),
\mbox{ in particular }
v_i \in L^\infty(\R^N).
\ee 

\subsection{Abstract setting for the boundary value problem \reff{BVP2}}
\label{abstract}
Let us introduce an abstract setting of the type of Corollary \ref{cor}
for the boundary value problem~\reff{BVP2}.  We take the
$\om_0$ from \reff{pom1}, 
$\om_1$ from  Theorem \ref{maxreg2} and $\om_2$ from \reff{vreg2}, and we fix an $\om$  with 
\bee
\label{omdef}
N-2<\om<\min\left\{\om_0,\om_1,\om_2\right\}.
\ee
The Banach spaces $U$ and $V$  and their norms are defined by
\bee
\label{UVdef2}
U:=W_0^{1,2,\om}(\Omega),\;
V:=W^{-1,2,\om}(\Omega),\;
\|\cdot\|_U:=\|\cdot\|_{1,2,\om},\;
\|\cdot\|_V:=\|\cdot\|_{-1,2,\om}.
\ee
The second norm $\|\cdot\|_\infty$ in $U$ of Corollary \ref{cor} is defined by \reff{normdef}.
Because of \reff{embed} the assumption \reff{weaker} of Corollary \ref{cor} is satisfied in this setting.
Further, the $C^1$-smooth operators $F_\ve:U \to V$ of Theorem \ref{ift} are defined by 
$$
F_\ve(u):=(A_\ve+B_\ve)u+C(u),
$$
where the 
nonlinear
operator $C:
L^\infty(\Omega)\to W^{-1,2,\om}(\Omega)$ is defined (similarly to \reff{Cdef1}) by
$$
\langle C(u),\vp\rangle_{1,2}:=
\int_\Omega\Big(c_i(x,u(x))\partial_{x_i}\vp(x)+d(x,u(x))\vp(x)\Big)dx
\mbox{ for all } \vp \in W_0^{1,2}(\Omega).
$$
Using this definition and  \reff{phidef1} we get that
\bee
\label{Cdef}
C(u)=\phi_{f(u),g(u)}
\mbox{ with } f_i(u):=c_i(\cdot,u(\cdot)),\;
g(u):=d(\cdot,u(\cdot))-u,
\ee
i.e. that $C$ is the superposition of the $C^1$-map
$u \in L^\infty(\Omega)\mapsto (f(u),g(u))\in L^\infty(\Omega)^N\times L^\infty(\Omega)$ (cf. \reff{diffass2}) and of the linear bounded map $(f,g) \in L^\infty(\Omega)^N\times L^\infty(\Omega)\mapsto \phi_{f,g}\in W^{-1,2,\om}(\Omega)$, i.e. $C$ is really $C^1$-smooth from
$L^\infty(\Omega)$ into $W^{-1,2,\om}(\Omega)$, and
\begin{eqnarray*}
&&\langle C'(u)v,\vp\rangle_{1,2}=\phi_{f'(u)v,g'(u)v}\nonumber\\
&&=\int_\Omega\Big(\partial_uc_i(x,u(x))\partial_{x_i}\vp(x)+\partial_ud(x,u(x))\vp(x)\Big)v(x)dx
\mbox{ for all } \vp \in W_0^{1,2}(\Omega).
\end{eqnarray*}
Hence, \reff{embed} implies that $C'(u)$ is compact from $W^{1,2,\om}(\Omega)$ into $W^{-1,2,\om}(\Omega)$ for any $u \in L^\infty(\Omega)$,
and, therefore, Theorem \ref{maxreg2} yields that the linear operators
$F'_\ve(u)=A_\ve+B_\ve+C'(u)$ 
are Fredholm of index zero from $U$ into $V$ for all $u \in U$.

Finally, we define the  approximate solution $u_0\in U$ of Corollary \ref{cor} to be the solution $u_0$ of the homogenized boundary value problem \reff{hombvp2}, which is given in Theorem \ref{main2}.
Remark that, similarly to \reff{unulleq1},
\bee
\label{unulleq2}
\hat Au_0+C(u_0)=0.
\ee

Because of Theorem \ref{maxreg2} 
we get that a function $u\in W_0^{1,2}(\Omega)\cap L^\infty(\Omega)$ is a weak solution to the boundary value problem 
\reff{BVP2} if and only if it and satisfies the operator equation $F_\ve(u)=0$. 
Hence, in order to prove Theorem \ref{main2}(i) by means of Corrolary \ref{cor} we have 
to verify assumption \reff{coerza}
and \reff{Fas1}
of  Corollary \ref{cor}
in the setting introduced above, i.e. that
there exists $\ve_0>0$ such that
\bee
\label{coerz2}
\inf\left\{
\|(A_\ve+B_\ve+C'(u_0))u\|_{-1,2,\om}:\;
\ve \in (0,\ve_0],\;
u \in W^{1,2,\om}(\Omega),\;
\|u\|_{1,2,\om}=1\right\}>0,
\ee
and that
\bee
\label{Cas2}
\sup_{\|v\|_{1,2,\om} \le 1}
\|(C'(u_0+u)-C'(u_0))v\|_{-1,2,\om}
\to 0
\mbox{ for } \|u\|_\infty \to 0,
\ee
and we have to construct a family $\uu \in W^{1,2,\om}(\Omega)$ such that there exists $\la>0$ such that
\bee
\label{newconve}
\|\uu-u_0\|_\infty+
\|(A_\ve+B_\ve)\uu+ C(\uu)\|_{-1,2,\om}=O(\ve^\la)
\mbox{ for } \ve \to 0.
\ee
This is what we are going to do below in the next three subsections.

\subsection{Verification of  \reff{coerz2}}
\label{sub3}
We proceed as in Subsection \ref{subcoerz1}. 
Suppose that \reff{coerz2} is not true. Then there exist sequences $\ve_1,\ve_2,\ldots>0$ and $u_1,u_2,\ldots \in W_0^{1,2,\om}(\Omega)$ such that 
\bee
\label{conve1}
\lim_{l\to \infty}\Big(\ve_l+\|(A_{\ve_l}+B_{\ve_l}+C'(u_0))u_l\|_{-1,2,\om}\Big)=0,
\ee
but
\bee
\label{normone1}
\|u_l\|_{1,2,\om}= 1 \mbox{ for all } l.
\ee
Because $W^{1,2,\om}(\Omega)$ is 
continuously embedded into $W^{1,2}(\Omega)$
 and because it is compactly embedded into  $L^\infty(\Omega)$, without loss of generality we may assume that there exist $u^*_1\in W^{1,2}(\Omega)$ and  $u^*_2\in L^\infty(\Omega)$ such that 
\bee
\label{infconv}
u_l\rightharpoonup u^*_1 \mbox{ for $l\to \infty$
weakly in }
W^{1,2}(\Omega)
\mbox{ and }
\lim_{l \to \infty}\|u_l-u^*_2\|_\infty=0.
\ee
For $\vp(x):=\mbox{sgn}(u^*_1(x)-u^*_2(x))$ follows that
\begin{eqnarray*}
0&=&\lim_{l \to \infty}\int_\Omega(u_l(x)-u^*_1(x))\vp(x)dx\\
&=&\lim_{l \to \infty}\int_\Omega(u_l(x)-u^*_2(x))\vp(x)dx+\int_\Omega|u^*_1(x)-u^*_2(x)|dx=\int_\Omega|u^*_1(x)-u^*_2(x)|dx,
\end{eqnarray*}
i.e.
$u^*_1=u^*_2=:u_*\in W^{1,2}(\Omega)\cap L^\infty(\Omega)$.
The functions $u_l$ are continuous and vanish on the boundary $\partial \Omega$,
and for $l\to \infty$ they converge uniformly on $\overline \Omega$ to the continous function $u_*$.
Therefore $u_*\in W_0^{1,2}(\Omega;\R^n)$.
Further, \reff{conve1} and
\reff{infconv}
imply that 
\bee
\label{conve}
\lim_{l\to \infty}\|(A_{\ve_l}+B_{\ve_l})u_l
+C'(u_0)u_*\|^*_{-1,2,\om}\to 0.
\ee 

The next lemma is the place where we use that the coefficients $b_{ij}$ are localized defects in the sense of assumption \reff{aass2}, i.e. that their supports are bounded:
\begin{lemma}
\label{Blemma2}
For any $0\le\rho<\sigma<N$ we have
$$
\sup\left\{\|B_{\ve}u\|_{-1,2,\rho}:\;
u \in W^{1,2,\sigma}(\Omega),\; \|u\|_{1,2,\sigma} 
\le 1\right\}=O(\ve^{(\sigma-\rho)/2})
\mbox{ for } \ve \to 0.
$$
\end{lemma}
{\bf Proof }
Because the supports of the functions $b_{ij}$ are bounded, there exists $R>0$ such that $b_{ij}(y)=0$ for $\|y\|>R$, i.e.
$
b_{ij}(\x)=0 \mbox{ for } x \notin \Omega_{0,\ve R}.
$ 
Here we use the notation \reff{Omxrdef}.
Take $\ve \in (0,1/R]$, $0\le\rho<\sigma<N$, $u \in W^{1,2,\sigma}(\Omega)$,
$r\in (0,1]$, $x \in \Omega$ and $\vp \in W_0^{1,2}(\Omega)$ with $\mbox{supp}\,\vp \subset \Omega_{x,r}$. Then the boundedness of the functions $b_{ij}$
yields that
\bee
\label{rest}
\left|\int_{\Omega}b_{ij}(\xi/\ve)\partial_{x_j}u(\xi)\partial_{x_i}\vp(\xi)d\xi\right|
\le \mbox{const} \, 
\left(\sum_{j=1}^N\int_{\Omega_{x,r}\cap\Omega_{0,\ve R}}\left|\partial_{x_j}u(\xi)\right|^2d\xi\right)^{1/2} \|\vp\|_{1,2},
\ee
where the constant  does not depend on $\ve$, $\rho$, $\sigma$,  $u$, $r$, $x$ and $\vp$. Moreover, because of $\partial_{x_j}u \in L^{2,\sigma}(\Omega)$ we have
$$
\int_{\Omega_{x,r}}\left|\partial_{x_j}u(\xi)\right|^2d\xi \le  r^\sigma\|u\|^2_{1,2,\sigma}
\mbox{ and }
\int_{\Omega_{0,\ve R}}\left|\partial_{x_j}u(\xi)\right|^2d\xi \le  (\ve R)^\sigma\|u\|^2_{1,2,\sigma}.
$$
Hence, \reff{rest} yields that
$$
r^{-\rho/2}
\left|\int_{\Omega}b_{ij}(\xi/\ve)\partial_{x_j}u(\xi)\partial_{x_i}\vp(\xi)d\xi\right|\le 
\mbox{const}\,
\min\left\{r^{(\sigma-\rho)/2},\ve^{\sigma/2} 
r^{-\rho/2}\right\}
\le 
\mbox{const}\,
\ve^{(\sigma-\rho)/2},
$$
where the constants do not depend on
$\ve$, $\rho$, $\sigma$,  $u$ and $r$.
Here we used that $r^{(\sigma-\rho)/2}\le 
\ve^{(\sigma-\rho)/2}$ for $r \le \ve$, and $\ve^{\sigma/2} 
r^{-\rho/2}\le 
\ve^{(\sigma-\rho)/2}$
for $r \ge \ve$.

Let us summarize: Take the
left-hand side if \reff{rest} and multiply it by $r^{-\rho/2}$ and take
the supremum over all
$u \in W^{1,2,\sigma}(\Omega)$
with $\|u\|_{1,2,\sigma}\le 1$,
$r\in (0,1]$, $x \in \Omega$ and $\vp \in W_0^{1,2}(\Omega)$ with $\mbox{supp}\,\vp \subset \Omega_{x,r}$
and $\|\vp\|_{1,2}\le 1$.
This way we get
$\sup\left\{\|B_{\ve}u\|_{-1,2,\rho}:\;
u \in W^{1,2,\sigma}(\Omega),\; \|u\|_{1,2,\sigma} 
\le 1\right\}$, and this is
$O(\ve^{(\sigma-\rho)/2})$
for $\ve \to 0$.
\qed\\

Now we use \reff{normone1} and
\reff{conve}
and  Lemma \ref{Blemma2} with 
$\rho=0$ and $\sigma=\om$ (cf. \reff{omdef}). This way it follows that 
\bee
\label{convea}
\lim_{l\to \infty}\|A_{\ve_l}
u_l+C'(u_0)u_*\|_{-1,2}=0.
\ee

In order to get the needed contradiction we use, as in Subsection \ref{subcoerz1}, Theorem \ref{Shentheorem1}.
Because of \reff{infconv}, \reff{convea} and Theorem \ref{Shentheorem1}
(with $f_*=C'(u_0)u_* \in W^{-1,2}(\Omega)$) 
it follows that
$
(\hat A+C'(u_0))u_*=0,
$
i.e. that
$u_*$ is a weak solution
to the linearized boundary value problem
\reff{linhombvp2}. Hence, by assumption of Theorem~\ref{main2}, we get that $u_*=0$.
Therefore
\reff{conve} 
implies that $\|(A_{\ve_l}
+B_{\ve_l})
u_l\|_{-1,2,\om}\to 0$ for $l \to \infty$.
But this contradicts to 
Theorem
\ref{maxreg2}
and \reff{normone1}.

\subsection{Verification of  \reff{Cas2}}
\label{sub3b}
We use the linear bounded 
map $(f,g)\in L^{\infty}(\Omega)^N
\times L^{\infty}(\Omega)
\mapsto \phi_{f,g}\in W^{-1,2,\om}(\Omega)$, which is defined by
\reff{phidef1}, and we use the representation \reff{Cdef} of the nonlinear operator $C$.
This way we get for $u,v \in W^{1,2,\om}(\Omega)$ that
$$
\big(C'(u_0+u)-C'(u_0)\big)v=
\phi_{(f'(u_0+u)-f'(u_0))v,(g'(u_0+u)-g'(u_0))v}.
$$
It follows that
\begin{eqnarray*}
&&\sup_{\|v\|_{1,2,\om} \le 1}\|(C'(u_0+u)-C'(u_0))v\|_{-1,2,\om}
\\
&&\le \mbox{const}
\sup_{\|v\|_\infty \le 1}
\|\phi_{(f'(u_0+u)-f'(u_0))v,(g'(u_0+u)-g'(u_0))v\|_{-1,2,\om}}\nonumber\\
&&\le\mbox{const}
\sup_{\|v\|_\infty \le 1}
\left(\sum_{i=1}^N\|(f_i'(u_0+u)-f'_i(u_0))v\|_\infty+\|(g'(u_0+u)-g'(u_0))v\|_\infty\right)\nonumber\\
&&\le\mbox{const}
\left(\sum_{i=1}^N\|f'_i(u_0+u)-f'_i(u_0)\|_\infty
+\|g'(u_0+u)-g'(u_0)\|_\infty\right),
\end{eqnarray*}
where the constants do not depend on  $u \in L^\infty(\Omega)$ with $\|u\|_\infty \le 1$. 
But
$$
\|f'_i(u_0+u)-f'_i(u_0)\|_\infty
=\esssup_{x \in \Omega}|\partial_uc_i(x,u_0(x)+u(x))-
\partial_uc_i(x,u_0(x))|
$$
tends to zero for $\|u\|_\infty \to 0$
because the functions $u\in \R \mapsto \partial_uc_i(\cdot,u)\in L^\infty(\Omega)$
are continuous and, hence, uniformly continuous on bounded sets
(cf. assumption \reff{diffass2}).
And similarly for 
$$
\|g'(u_0+u)-g'(u_0)\|_\infty
=\esssup_{x \in \Omega}|\partial_ud(x,u_0(x)+u(x))-
\partial_ud(x,u_0(x))|.
$$
Hence, assumption \reff{Fas1} of Corrolary \ref{cor} is satisfied in the setting introduced in Subsection \ref{abstract}.

\subsection{Construction of approximate solutions with \reff{newconve}}
\label{sub1}
In this subsection we proceed as in Subsection \ref{sublim1}.
We define  a family $\uu\in  W^{1,2,\om}(\Omega)$, which satisfies 
\reff{newconve}, similarly to \reff{barudef1} as follows:
\bee
\label{barudef2}
\uu(x):=
u_0(x)+\ve  \eta_\ve(x)\partial_{x_k}u_0(x) v_k(\x).
\ee
Here $\eta_\ve$ is a cut-off function as in \reff{etaest},
and
$v_k\in W^{1,2}_{\rm loc}(\R^N)$ are the the $\Z^N$-periodic correctors, which are defined by the cell problems~\reff{cell2}. 
Because of \reff{pom1}, \reff{vreg2}
and \reff{barudef2} it  follows that
\bee
\label{pnull}
\uu \in W^{1,2,\om_*}(\Omega).
\ee
In particular, we have that $\uu \in W^{1,2,\om}(\Omega)$ with $\om$ defined in \reff{omdef}, as needed.
Further, from  \reff{pom1} and
\reff{vreg2}
follows
\bee
\label{barest}
\|\uu-u_0\|_\infty \le \mbox{const }\ve\|u_0\|_{1,2},
\ee
where the constant does not depend on $\ve$.
Hence, in order to prove \reff{newconve}
it remains to show that there exists $\la>0$ such that
$\|(A_\ve+B_\ve)\uu+C(u_0)\|_{-1,2,\om}=O(\ve^\la)$ for $\ve \to 0$ or, what is that same
(because of \reff{unulleq2}), that
$$
\|(A_\ve+B_\ve)\uu-\hat Au_0\|_{-1,2,\om}=O(\ve^\la) \mbox{ for } \ve \to 0.
$$
But this is the assertion of
the two lemmas below.
\begin{lemma}
\label{BLemma2}
There exists $\la>0$ such that
$\|B_\ve\uu\|_{-1,2,\om}=O(\ve^\la)$ for $\ve \to 0$.
\end{lemma}
{\bf Proof }
We proceed as in the proof of Lemma \ref{BLemma1}.
Because of Lemma \ref{Blemma2}  with $\rho=\om$ and $\sigma=\om_*$ 
(cf. \reff{omdef} and \reff{pnull}) we have that
$$
\|B_\ve\uu\|_{-1,2,\om} =O\Big(\ve^{(\om_*-\om)/2}\|\uu\|_{1,2,\om_*}\Big)
\mbox{ for } \ve \to 0.
$$
Hence, it sufficies to show
that there exists $\la>0$ such
that
$\ve^{(\om_*-\om)/2}\|\uu\|_{1,2,\om_*}=O(\ve^\la)$ for $\ve \to 0$ or, what is the same, that
there exists $\la>0$ such that
$\ve^{(\om_*-\om)/2}
\|\partial_{x_i}\uu-\partial_{x_i}u_0)\|_{2,\om_*}=O(\ve^\la)$ for $\ve \to 0$.
Using 
\begin{eqnarray*}
&&\partial_{x_j}\uu(x)-\partial_{x_j}u_0(x)\\
&&=\eta_\ve(x) 
\partial_{x_k}u_0(x) \partial_{y_j}v_k(\x)+\ve\Big(\partial_{x_j}\eta_\ve(x)\partial_{x_k}u_0(x) v_k(\x)+
\eta_\ve(x)
\partial_{x_j}\partial_{x_k}u_0(x) v_k(\x)\Big),
\end{eqnarray*}
we get finally that it remains to prove that there exists $\la>0$ such that
\begin{eqnarray}
\label{este1}
&&\ve^{\om_*-\om}
\sup_{x \in \Omega,r \in (0,1]} r^{-\om_*}
\int_{\Omega_{x,r}} \left|
\eta_\ve(\xi)\partial_{x_k}u_0(\xi) \partial_{y_j}v_k(\xi/\ve)\right|^2d\xi=O(\ve^\la) \mbox{ for } \ve \to 0,\\
\label{este2}
&&\ve^{2+\om_*-\om}
\sup_{x \in \Omega,r \in (0,1]} r^{-\om_*}
\int_{\Omega_{x,r}} \left|\partial_{x_j}\eta_\ve(\xi)\partial_{x_k}u_0(\xi)v_k(\xi/\ve)\right|^{2}d\xi=O(\ve^\la) \mbox{ for } \ve \to 0,
\\
\label{este3}
&&\ve^{2+\om_*-\om}
\sup_{x \in \Omega,r \in (0,1]} r^{-\om_*}
\int_{\Omega_{x,r}} \left|\eta_\ve(\xi)\partial_{x_j}\partial_{x_k}u_0(\xi)v_k(\xi/\ve)\right|^{2}d\xi=O(\ve^\la) \mbox{ for } \ve \to 0.
\end{eqnarray}

In a first step, let us prove the assertion \reff{este1}. Because of 
\reff{pom} we have for all $\ve>0$ and all $\xi \in \R^N$ that
\bee
\label{1est}
|\eta_\ve(\xi)\partial_{x_k}u_0(\xi) \partial_{y_j}v_k(\xi/\ve)|
\le \mbox{const } 
|\partial_{y_j}v_k(\xi/\ve)|,
\ee
where the constant does not depend on $\ve$ and $\xi$. 
Take $x \in \Omega$ and $r>0$.
Because the functions $\partial_{y_j}v_k$ are $\Z^N$-periodic and belong to $L^{2,\om_*}_{\rm loc}(\R^N)$ and because of \reff{cu} we have
$$
\int_{\|\xi-x\|< r} |\partial_{y_j}v_k(\xi/\ve)|^{2}d\xi\le \ve^N
\int_{\|y\|< r/\ve} |\partial_{y_j}v_k(x/\ve+y)|^{2}dy
\le  \mbox{const}\,
\left\{
\begin{array}{ccc}
\ve^N(r/\ve)^{\om_*} &\mbox{if}& 
r\le\ve,\\
(r+\ve)^N &\mbox{if}& r\ge\ve
\end{array}
\right.
$$
and, hence,
\bee
\label{rve}
r^{-\om_*}\int_{\|\xi-x\|< r} |\partial_{y_j}v_k(\xi/\ve)|^{2}d\xi
\le \mbox{const}\,\ve^{N-\om_*},
\ee
where the constants do not depend on $x$, $r$ and  $\ve$. 
This way we get from \reff{1est}
$$
\ve^{\om_*-\om}
\sup_{x \in \Omega,r \in (0,1]} r^{-\om_*}
\int_{\Omega_{x,r}} |\eta_\ve(x)\partial_{x_k}u_0(\xi) \partial_{y_j}v_k(\xi/\ve)|^2d\xi
\le \mbox{const} \,\ve^{N-\om},
$$
where the constant does not depend $\ve$.
Hence, \reff{este1} is verified.

In a second step, let us prove the assertion \reff{este2}. Because of \reff{etaest}, \reff{pom} and \reff{vreg2} we have
$$
\ve^{2+\om_*-\om} r^{-\om_*}
\int_{\Omega_{x,r}} \left|\partial_{x_j}\eta_\ve(\xi)\partial_{x_k}u_0(\xi)v_k(\xi/\ve)\right|^{2}d\xi
\le \mbox{const }\ve^{\om_*-\om} r^{-\om_*}
|\Omega_{x,r} \cap \Omega_\ve|,
$$
where the constant does not depend on $x$, $r$ and  $\ve$. 
But \reff{Omest} yields that
\bee
\label{rve1}
r^{-\om_*}
|\Omega_{x,r}\cap \Omega_\ve|\le 
\left\{
\begin{array}{cl}
r^{-\om_*}|\Omega_{x,r}|\le r^{N-\om_*}\le \ve^{1-\om_*/N} & \mbox{ if } r \le \ve^{1/N},\\
\ve^{-\om_*/N}|\Omega_{\ve}|\le \ve^{1-\om_*/N} & \mbox{ if } r \ge \ve^{1/N}.
\end{array}
\right.
\ee
Therefore
$$
\ve^{2+\om_*-\om} 
\sup_{x \in \Omega,r \in (0,1]}
r^{-\om_*}
\int_{\Omega_{x,r}} \left|\partial_{x_j}\eta_\ve(\xi)\partial_{x_k}u_0(\xi)v_k(\xi/\ve)\right|^{2}d\xi
\le\mbox{const }\ve^{\om_*-\om+1-\om_*/N},
$$
where the constant does not depend on $\ve$.

And finally, let us prove \reff{este3}. Because of \reff{pom1} and of the boundedness of $v_k$ we have for $x \in \Omega$ and $r>0$ that
$$
\int_{\Omega_{x,r}}|\eta_\ve(x)\partial_{x_j}\partial_{x_k}u_0(\xi) v_k(\xi/\ve)|^2d\xi
\le \mbox{const } r^{\om_*},
$$
where the constants do not depend on $x$, $r$ and $\ve$. It follows that
$$
\ve^{2+\om_*-\om}
\sup_{x \in \Omega,r \in (0,1]} r^{-\om_*}
\int_{\Omega_{x,r}} |\eta_\ve(x)\partial_{x_j}\partial_{x_k}u_0
(\xi) v_k(\xi/\ve)|^2d\xi
\le\mbox{const }  \ve^{2+\om_*-\om},
$$
where the constant does not depend on $\ve$.
\qed

\begin{lemma}
\label{Klemma2}
We have
$\|A_\ve\uu-\hat Au_0\|
_{-1,2,\om}=O(\ve^\la)$
for $\ve \to 0$
with $\la:=\frac{1}{2}\left(1-\frac{\om}{N}\right)$.
\end{lemma}
{\bf Proof }
We proceed as in the proof of Lemma \ref{Klemma}.

For $i,j,k=1,\ldots,N$ 
we define
$\Z^N$-periodic functions $f_{ij}\in L_{\rm loc}^{2}(\R^N)$ and $g_{ij}\in W_{\rm loc}^{2,2}(\R^N)$
and 
$h_{ijk}\in W_{\rm loc}^{1,2}(\R^N)$ 
by
\begin{eqnarray*}
&&f_{ij}(y):=a_{ij}(y)
+a_{ik}(y)
\partial_{y_k}v_j(y)-\hat a_{ij},\\
&&
\Delta g_{ij}(y)=f_{ij}(y),\; \int_{(0,1)^N}g_{ij}(y)dy=0,\\
&&
h_{ijk}(y):=
\partial_{y_i}g_{jk}(y)
-\partial_{y_j}g_{ik}(y).
\end{eqnarray*}
Because of \reff{vreg2}
we have that
$f_{ij}\in L_{\rm loc}^{2,\om}(\R^N)$, and, hence, 
$\partial_{y_k}g_{ij}, h_{ijk} \in W_{\rm loc}^{1,2,\om}(\R^N)$, in particular
\bee
\label{hbound2}
h_{ijk} \in L^{\infty}(\R^N)
\ee
(cf. also \cite[Remark 2.1]{Kenig}).
As is \reff{threeint1} we get 
for any $\vp \in W_0^{1,2}(\Omega)$ that
\begin{eqnarray}
&&
\langle A_\ve\uu-\hat Au_0,\vp\rangle_{1,2}
\nonumber\\
&&=\ve\int_\Omega\left(-h_{ijk}^{\al \beta}(\x)+
a^{\al\gamma}_{ik}(\x) v_j^{\gamma\beta}(\x)
\right)
\partial_{x_k}\Big(\eta_\ve(x)
\partial_{x_j}u_0^\beta(x)\Big)
\partial_{x_i}\vp^\al(x)dx\nonumber\\
&&\;\;\;\;\;\;+\int_\Omega
\left(a^{\al \beta}_{ij}(\x)
-\hat a^{\al \beta}_{ij}\right)\left(1-\eta_\ve(x)\right)
\partial_{x_j}u_0^\beta(x)
\partial_{x_i}\vp^\al(x)dx.
\label{threeint}
\end{eqnarray}
We
estimate
$\| A_\ve\uu-\hat Au_0\|_{-1,2,\om}$
by means of \reff{threeint}.
We have to estimate the supremum of $r^{-\om/2}\langle A_\ve\uu-\hat Au_0,\vp\rangle_{1,2}$ over all $x \in \Omega$, $r\in(0,1]$,
 and $\vp \in W_0^{1,2}(\Omega)$ with $\mbox{supp} \,\vp \subset \Omega_{x,r}$ and $\|\vp\|_{1,2}\le 1 $.
 
Take $\ve>0$,  $x \in \Omega$, $r\in(0,1]$
and $\vp \in W_0^{1,2}(\Omega)$ with $
\mbox{supp} \,\vp \subset \Omega_{x,r}$
and $\|\vp\|_{1,2}\le 1$. 
In what follows we will use 
the fact that the functions $a_{ij}$, $v_k$ and $h_{ijk}$ are bounded
(cf. \reff{aass2}, \reff{vreg2} and \reff{hbound2}), and all constants will be independent on $\ve$, $x$, $r$
and $\vp$.
 
We split the right-hand side of  \reff{threeint}, multiplied by $r^{-\om/2}$, into three terms. The first term can be estimated as follows:
\begin{eqnarray*}
&&r^{-\om/2}\left|
\int_\Omega
\left(a_{ij}(\xi/\ve)
-\hat a_{ij}\right)\left(1-\eta_\ve(x)\right)
\partial_{x_j}u_0(\xi)
\partial_{x_i}\vp(\xi)d\xi\right|\\
&&\le \mbox{const }r^{-\om/2}
|\Omega_{x,r}\cap \Omega_\ve|^{1/2}
\le \mbox{const } \ve^\la \mbox{ with } \la=\frac{1}{2}\left(1-\frac{\om}{N}\right).
\end{eqnarray*}
Here we used \reff{rve1}. 
By means of \reff{pom1} the 
second term can be estimated as 
\begin{eqnarray*}
&& \ve r^{-\om/2}\left|\int_{\Omega}
\Big(a_{ik}(y/\ve) v_j(\xi/\ve)-h_{ijk}(\xi/\ve)\Big)\eta_\ve(x)
\partial_{x_k}\partial_{x_j}u_0
(\xi)\partial_{x_i}\vp(\xi)d\xi\right|
\nonumber\\
&&
\le\mbox{const } 
\ve r^{-\om/2}
\left(\sum_{j,k=1}^N\int_{\Omega_{x,r}}|[\partial_{x_k}\partial_{x_j}
u_0(\xi)|^2d\xi\right)^{1/2}
\le\mbox{const }\ve.
\end{eqnarray*}
And the 
third term can be estimated, using \reff{etaest} and \reff{rve1}, as 
\begin{eqnarray*}
&&\ve r^{-\om/2}\left|\int_{\Omega}
\Big(a_{ik}(y/\ve) v_j(\xi/\ve)-h_{ijk}(\xi/\ve)\Big)\partial_{x_j}\eta_\ve(x)
\partial_{x_k}u_0
(\xi)\partial_{x_i}\vp(\xi)d\xi\right|\\
&& \le \mbox{const }
\ve r^{-\om/2}\left(\sum_{k=1}^n\int_{\Omega_{x,r}\cap \Omega_\ve}
|\partial_{x_j}\eta_\ve(x)
\partial_{x_k}u_0
(\xi)|^2d\xi\right)^{1/2}\\
&&\le \mbox{const }
r^{-\om/2}|\Omega_{x,r}\cap \Omega_\ve|^{1/2}\le \mbox{const }
\ve^\la \mbox{ with } \la=\frac{1}{2}\left(1-\frac{\om}{N}\right).
\end{eqnarray*}
\qed\\

\section{Proof of Theorem  \ref{main2}(ii)}
\label{sec:proof2ii}
\setcounter{equation}{0}
\setcounter{theorem}{0}

Because of \reff{omdef} the exponent $\om>N-2$ may be choosen arbitrarily close to $N-2$. Therefore
$\la=\frac{1}{2}\left(1-\frac{\om}{N}\right)$
can be made to be arbitrarily close to $1/N$. Hence, \reff{barest} and Lemma \ref{Klemma2}
yield the assertion of Theorem  \ref{main2}(ii).

\section*{Acknowledgments}

The author gratefully acknowledges that
his paper is the result of 
a long-standing mathematical cooperation and friendship with Jens A. Griepentrog.


\begin{thebibliography}{troi}


\bibitem{Barbu}
V. Barbu, Partial Differential Equations and Boundary Value Problems. Mathematics and its  Applications vol. {\bf 441}, Kluwer Academic Publishers, 1998. 


\bibitem{BenF}
A. Bensoussan, J. Frehse, Regularity Results for Nonlinear Elliptic Systems and Applications. Applied Mathematical Sciences vol. {\bf 151}, Springer, 2002. 


\bibitem{Ben}
A. Bensoussan, J.L. Lions, G. Papanicolaou, Asymptotic Analysis for Periodic Structures. Studies in Mathematics and its  Applications vol. {\bf 3}, North-Holland, 1978. 

\bibitem{BJL}
X. Blanc, M. Josien, C. Le Bris,
Precise approximations in elliptic homogenization beyond periodic setting. Asympttot. Anal. {\bf 116} (2020), 93--137.



\bibitem{Blanc}
X. Blanc, C. Le Bris,
Homogenization Theory for Multiscale Problems. An Introduction. 
Modeling, Simulation and Applications vol. {\bf 21}, Springer, 2023.

\bibitem{BDL} A. Braides, G. Dal Maso, C. Le Bris, A closure theorem for $\Gamma$-convergence and H-convergence with applications to non-periodic homogenization.
arXiv:2402.19031.

\bibitem{Breden} M. Breden, R. Castelli, Existence and instability of steady states
for a triangular cross-diffusion system: a computer-assisted proof.
J. Differ. Equations {\bf 264} (2018), 6418--6458.

\bibitem{Bun} R. Bunoiu, R. Precup, Localization and multiplicity in the homogenization of nonlinear problems.
Adv. Nonlinear Anal. {\bf 9} (2020), 292--304.




\bibitem{Butetc} V.F. Butuzov, N.N. Nefedov, O.E. Omel'chenko, L. Recke,
Time-periodic boundary layer solutions to singularly perturbed parabolic problems.
J. Differ. Equations {\bf 262} (2017), 4823--4862.

\bibitem{But2022}  
V.F. Butuzov, N.N. Nefedov,  O.E. Omel'chenko, L. Recke, 
Boundary layer solutions to singularly perturbed quasilinear systems.
Discrete Cont. Dyn. Syst., Series B {\bf 27}
(2022), 4255--4283.


\bibitem{Fiedler}  
V.F. Butuzov, N.N. Nefedov,  O.E. Omel'chenko, L. Recke, K.R. Schneider, 
An implicit function theorem and applications to nonsmooth boundary layers.
In: {\it Patterns of Dynamics}, ed. by 
P. Gurevich, J. Hell, B. Sandstede, A. Scheel, Springer Proc. in Mathematics \& Statistics vol. {\bf 205}, Springer,  
2017, 111--127.


\bibitem{Ca} 
M. Cadiot, J.-P. Lessard,  J.-Ch. Nave, 
Stationary non-radial localized patterns in the planar Swift-Hohenberg PDE: constructive proofs of existence.
J. Differ. Equations {\bf 414} (2025), 555--608.






\bibitem{CP}
L. Caffarelli, I. Peral,
On $W^{1,p}$ estimates for elliptic equations in divergence form.
Comm. Pure Appl. Math. {\bf 51} (1998) 1--21.


\bibitem{Che}
G.A. Chechkin, A.L. Piatnitski, A.S. Shamaev, Homogenization. Methods and Applications.
Translations of Mathematical Monographs vol. {\bf 234}, AMS, Providence 2007. 



\bibitem{Chen}
Ya-Zhe Chen, Lan-Cheng Wu, Second-Order Elliptic Equations and Elliptic Systems.
Translations of Mathematical Monographs {\bf 174}, Providence, 1998 (Chinese Original: Beijing, 1991).



\bibitem{Ci}
D. Cioranescu, P. Donato, An Introduction to Homogenization. Oxford Lecture Series in Mathematics and its Applications vol. {\bf 17}, Oxford University Press, 1999. 






\bibitem{Ga}
T. Gallouet, A. Monier, On the regularity of solutions of elliptic equations. 
Rendiconti di Matematica {\bf VII 19}
(1999), 471--488. 




\bibitem{Gia}
M. Giaquinta, Introduction to Regularity Theory for Nonlinear Elliptc Systems. Lecture Notes in Mathematics, ETH Z\"urich, Birkh\"auser, 1993.

\bibitem{Giusti}
M. Giusti, Direct Methods in the Calculus of Variations. World Scientific, 2003.



\bibitem{Griep} 
J. A. Griepentrog,
Linear elliptic boundary value problems with non-smooth data: Campanato spaces of functionals. 
Math. Nachr. {\bf 243} (2002), 19--42. 





\bibitem{GR}
J.A. Griepentrog, L. Recke,
Linear elliptic boundary value problems with non-smooth data: Normal solvability in Sobolev-Campanato spaces. 
Math. Nachr. {\bf 225} (2001), 39--74. 

\bibitem{GREvol}
J.A. Griepentrog, L. Recke,
Local existence, uniqueness and smooth dependence for nonsmooth quasilinear parabolic problems. 
J. Evol. Equ. {\bf 10} (2010), 341-375. 


\bibitem{G} K. Gr\"oger, A $W^{1,p}$-estimate for solutions to mixed boundary value problems for second-order elliptic differential equations. Math. Ann. 
{\bf 283} (1989), 679--687.


\bibitem{GrR} K. Gr\"oger, L. Recke,
Applications of differential calculus to quasilinear elliptic boundary value problems with non-smooth data. 
NoDEA, Nonlinear Differ. Equ. Appl. {\bf 13}
(2006), 263--285.




\bibitem{Katz}
S.G. Katz, H.R. Parks, The Implicit Function Theorem. History, Theory, and Applications.
Birkh\"auser, 2002.

\bibitem{Kenig}
C. Kenig, Fanghua Lin,  Zhongwei Shen, 
Periodic homogenization of Green and Neumann functions. 
Commun. Pure Appl. Math. {\bf 67} (2014), 1219-1262.





\bibitem{Kufner}
A. Kufner, O. John, S. Fucik, Function Spaces.
Academia, Prague, 1977.



\bibitem{Lanza1} M. Lanza de Cristoforis, P. Musolino, Two-parameter homogenization for a nonlinear periodic Robin problem for a Poisson equation: a functional analytic approach. Rev. Mat. Complut. {\bf 31} (2018), 63--110.

\bibitem{Lanza2} M. Lanza de Cristoforis, P. Musolino, Asymptotic behaviour of the energy integral of a two-parameter homogenization
problem with nonlinear periodic Robin boundary conditions. Proc. Edinb. Math. Soc. II. Ser. {\bf 62} (2019), 985--1016.


\bibitem{Magnus1} R.J. Magnus, The implicit function theorem and multi-bump solutions of periodic partial differential equations. Proc.
Royal Soc.  Edinb.  {\bf 136A} (2006), 559--583.


\bibitem{Magnus2} R.J. Magnus, A scaling approach to bumps and multi-bumps for nonlinear partial differential 
equations. Proc.
Royal Soc.  Edinb.  {\bf 136A} (2006), 585--614.



\bibitem{ME1975}
N.G. Meyers, An $L^p$-estimate for the gradient of solutions of second-order elliptic divergence equations. Ann. Scuola Norm. Sup. Pisa {\bf 17} (1963), 189--206.


\bibitem{M}
N.G. Meyers, A. Elcrat,
Some results on regularity for solutions of non-linear elliptic systems and quasi-regular functions.
Duke Math. J.  {\bf 42} (1975), 121--136.







\bibitem{NURS}  
N.N. Nefedov, A.O. Orlov, L. Recke, K.R. Schneider, Nonsmooth regular perturbations of singularly perturbed problems. 
J. Differ. Equations {\bf 375} (2023), 206--236.


\bibitem{I} N.N. Nefedov, L.~Recke, A common approach to singular perturbation and homogenization I: Quasilinear ODE systems. arXiv:2309.15611.


\bibitem{II} N.N. Nefedov, L.~Recke, A common approach to singular perturbation and homogenization II: Semilinear elliptic PDE systems. 
J. Math. Anal. Appl. {\bf 545} (2025), Article ID 129099.



\bibitem{OR2009}
O.E Omel'chenko, L. Recke, 
Boundary layer solutions to singularly perturbed problems via the implicit function theorem. 
Asymptotic Anal. {\bf 62} (2009), 207-225.



\bibitem{OmelchenkoRecke2015} O.E.~Omel'chenko, L.~Recke,
Existence, local uniqueness and asymptotic approximation of spike solutions
to singularly perturbed elliptic problems. Hiroshima Math. J. {\bf 45} (2015), 35--89.

\bibitem{PRS}
D.K. Palagachev, L. Recke, L.G. Softova, 
Applications of the differential calculus to nonlinear elliptic operators with discontinuous coefficients. 
Math. Ann. {\bf 336} (2006), 617-637. 







\bibitem{Recke1995} L.~Recke,
 Applications of the implicit function theorem to quasilinear elliptic boundary value problems with non-smooth data. 
Commun. Partial Differ. Equations {\bf 20} (1995), 1457-1479. 




\bibitem{Recke2022} L.~Recke, Use of very weak approximate
boundary layer solutions to spatially nonsmooth
singularly perturbed problems.
J. Math. Anal. Appl. {\bf 506} (2022), Article ID 125552.

\bibitem{ReckeNonper} L.~Recke, Nonlinear non-periodic homogenization:
Existence, local uniqueness and estimates.
arXiv:2408.06705.


\bibitem{ReckeOmelchenko2008} L.~Recke, O.E.~Omel'chenko,
Boundary layer solutions to problems
with infinite dimensional singular and regular perturbations.
J. Differ. Equations {\bf 245} (2008), 3806--3822.







\bibitem{Riva}
M.D. Riva, R. Molinarolo, P. Musolino, 
Local uniqueness of the solutions for a singularly perturbed nonlinear nonautonomous transmission problem. 
Nonlinear Anal., Theory Methods Appl., Ser. A {\bf 191} (2020), Article ID 111645. 




\bibitem{Sa}
Y. Sawano, G. Di Fazio, D.I. Hakim, Morrey Spaces. Introduction and Applications to Integral Operators and PDEs, Volumes I and II. Monographs and Research Notes in Math., CPC Press, 2020.






\bibitem{Shen}
Zongwei Shen, Periodic Homogenization of Elliptic Systems. Operator Theory: Advances and Applications vol. {\bf 269}, Birkh\"auser, 2018.


\bibitem{Tartar}
L. Tartar,
The General Theory of Homogenization. A Personalized Introduction. 
Lecture Notes of the Unione Matematica Italiana vol. {\bf7}, Springer, 2009. 











\bibitem{Troi}
G.M. Troianiello, Elliptic Differential Equations and Obstacle Problems. Plenum Press, 1987.













\end{thebibliography}
\end{document}